\documentclass[a4paper,12pt]{article}
\usepackage{bm}
\usepackage{amsmath}
\usepackage{amssymb}
\usepackage[dvips]{graphicx}
\usepackage{pdflscape}
\usepackage{multirow}
\usepackage[sort]{natbib}
\usepackage{soul}
\usepackage{xr}
\usepackage{amsmath}
\usepackage{graphicx}
\usepackage{enumerate}
\usepackage{natbib}
\usepackage[hyphens]{url} 
\usepackage[utf8]{inputenc}
\usepackage{mathtools}
\usepackage{amsbsy}
\usepackage{amssymb}
\usepackage{color}
\usepackage{ulem}
\usepackage{placeins}
\usepackage{booktabs}
\usepackage{multirow}
\usepackage{hyperref}
 \usepackage{booktabs,subcaption,amsfonts,dcolumn}
\usepackage{xcolor,colortbl}
\usepackage{subcaption}
\usepackage{verbatim}
\usepackage{chngcntr}
\usepackage{apptools}

\usepackage{amsthm}
\newtheorem{thm}{Theorem}

\newtheorem{defin}{Definition}
\newtheorem{eg}{Example}

\definecolor{Gray}{gray}{0.85}
\definecolor{LightCyan}{rgb}{0.88,1,1}
 
\newcolumntype{b}{>{\columncolor{Gray}}c}
\newcolumntype{a}{>{\columncolor{red}}c}
\newcolumntype{g}{>{\columncolor{green}}c}

\title{Families of discrete circular distributions with some novel applications \footnote{An early version of this work was presented at the Leeds Annual Statistics Research (LASR) Workshop 2019.}}
\author{ Kanti V. Mardia\thanks{%
University of Leeds and University of Oxford.} \and Karthik Sriram\thanks{%
Indian Institute of Management Ahmedabad, India.}}

\begin{document}
\date{} \maketitle
\begin{abstract}
We give a unified treatment of constructing families of circular discrete distributions.  Some of these families are deduced from established distributions such as von Mises and wrapped Cauchy. Some others are derived directly such as a flexible family based on trigonometric sums and the circular location family.  Results interrelating these families are discussed.  These distributions have been motivated by  two examples of discrete circular data: casino roulette spins and smart health acrophase monitoring, and   these data are analyzed using our proposed models.  We discuss how using continuous circular  models for circular discrete data can be misleading.
\end{abstract}

\section{Introduction}
\label{S:intro}
The subject of Directional Statistics has grown tremendously, especially since the 1980's, with advances in  ``Statistics on Manifolds" leading to new distributions on the hyper-sphere, torus, Stiefel manifold, Grassmann manifold and so on. The progress in this area  can be seen through several books published since then:  \citet{fisheretal1987}, \citet{fisher1993}, \citet{mardiajupp2000},   \citet{jammala_seng2001}, \citet{ley_verdebout2017} and \citet{ley_verdebout2018}. There has been a  recent special issue of Sankhya   edited by \citet{spissue2019}.  Further, \citet{pewsey_eduardo2020} have given a comprehensive survey of directional statistics and in the discussion to the paper, \citet{mardia2021} has given a brief history of the subject. However, there is limited development of circular discrete models.  There are good choices for continuous models for circular data, but  there has been  a dearth of  models for discrete data.  In this paper, we give the first unified treatment of constructing families of circular discrete distributions and present examples of circular data  that are observed directly as discrete rather than created by  grouping continuous data. The two data that motivated our paper are:
 
~\\
(i) {\bf Roulette wheel data:}  A typical European roulette  has $37$ discrete outcomes, viz. $\{0,1,2,\ldots, 36\}$. If the outcome  $0$ is mapped to $0$ radians, then the outcomes get mapped to a regular support of $37$ points on the circle given by
\begin{eqnarray}
\left\{ \frac{2\pi r}{37}, ~r\in \{0,1,2,\ldots, 36\}\right\}.\label{eq:supp_rou}
 \end{eqnarray}
 In Section \ref{S:examples_roulette}, we consider data sequences obtained from spins of four different European roulette wheels, one from an online roulette simulator and three from two different casino industries.  It should be pointed out that Karl Pearson, in the early 1890's, acquired   roulette spins data from Monte Carlo  to examine the question of whether the roulette wheel was unbiased (see, \citealt[p.60]{plackett1983}), and indeed  his paper of 1897 has the apt title ``The scientific aspect of Monte Carlo roulette" (\citealt{pearson1897a}).  Few authors have  considered this problem but have used  linearized methods, beginning with Karl Pearson and subsequently some others,  e.g.  \citet{ethier1982}, \citet{spencer2009}. Surprisingly, for this important application from the gaming industry, there has been little attention to inference that uses explicitly the circularity of the wheel.
 
 ~\\
(ii) {\bf Acrophase data:} In non-invasive smart health monitoring,  parameters such as Systolic Blood Pressure (SBP) are  recorded, by ambulatory devices, at predetermined discrete time points repeated each day. ``Acrophase" is defined as the time point at which the maximum SBP reading is recorded on a given day. Typically, acrophase data is  extracted from SBP measurements at each half hour during daytime (8 am to 8 pm) and  each hour during nighttime (8 pm to 8 am).  If we map 8 am to 0 radians and 8 pm to $\pi$ radians, the acrophase times get mapped to an irregular support of $36$ points on the circle given by
\begin{eqnarray}
\left\{\frac{2\pi r}{48}, r =0,1,2,\ldots,24\right\} \bigcup \left\{\frac{2\pi r}{48}, r =26, 28,\ldots, 44, 46\right\},  \label{eq:supp_acro}
 \end{eqnarray}
where the first set in the union corresponds to 25 half-hourly  points during daytime and the second set corresponds to the 11 one-hourly points during nighttime. 

In any application with  discrete data, one usually takes into account the  discrete nature of the underlying population. Our overall  recommendation is ``if one has discrete circular  data then one should start with a discrete circular model".  Also, the ``loss" due to use of a continuous model for circular discrete data can only be assessed after appropriate  discrete modeling, which  serves as a benchmark. As an illustration, for the acrophase data, we see in Section \ref{S:acrophase} that ignoring the underlying discreteness can lead to biased estimates.  Of course, this issue of discrete versus continuous distributions is a general problem, which is well known and has been dealt successfully in Linear Statistics but we need to treat this problem here as a model misspecification  (see, Section \ref{S:comparisons}).

In this paper, we give four methods to construct families of discrete distributions on the circle along with some basic results interrelating the methods.  We apply these models to analyze the aforementioned examples of  discrete data and also to provide insights  based on comparisons among discrete as well as (approximate) continuous models for discrete data. Our  methods to construct the probability distributions can be briefly described as follows:
 \begin{itemize}
 \item[(i)] Maximum entropy method: We start with a given set of  moment conditions for the discrete distribution on the circle. We then determine the discrete probability distribution with the maximum Shannon entropy among those satisfying the moment constraints. 
 \item[(ii)] Centered wrapping method: We start with a given discrete distribution on the line, and wrap it on the circle to obtain a discrete distribution on the circle. 
 \item[(iii)] Marginalized method: We start with a continuous  distribution on the circle, which we refer to  as the ``parent", and then obtain a discrete distribution on the circle by integrating the probability density function (pdf) on pre-determined arcs on the circle.
 \item[(iv)] Conditionalized method: We start with a continuous  distribution on the circle (parent), and then obtain the discrete distribution on the circle by restricting and normalizing the pdf to a pre-determined  lattice on the circle.
 \end{itemize}
In particular, we derive circular discrete distributions from general continuous location families and a family based on trigonometric sums. Key special cases include discrete families deduced from two established continuous distributions, viz. von Mises and wrapped Cauchy. These two distributions are commonly selected for circular data depending on whether the unimodal data has a long tail (wrapped Cauchy) or not (von Mises), which we now describe.

The direction of a unit random vector in two dimensions can be represented by an angle $\Theta$. On the circle, the von Mises distribution for $\Theta$ (see, for example, \citealt[p.36]{mardiajupp2000}) plays the same role as  the  normal distribution on the line. It belongs to the exponential family with two analogous parameters. Its pdf is given by
\begin{equation}
\label{eq:von_Mises1}
f_v(\theta\vert \kappa, \mu) = \frac{1}{2\pi I_0(\kappa)} e^{\kappa \cos(\theta-\mu)}, \theta\in [0, 2\pi), \mu\in [0, 2\pi), \kappa\geq 0,
\end{equation} 
where $\mu$ is the mean direction and $\kappa$ is the concentration (precision) parameter. The normalization constant $I_0(\kappa)$ is the modified Bessel function of order $0$:
$$ I_0(\kappa)= \sum_{r=0}^\infty \frac{\kappa^{2r}}{(r !)^2}.$$ 
 For large $\kappa$,  $\Theta$  is approximately normal  with mean $\mu$ and  variance 2/$\kappa$, and for $\kappa=0$, $\Theta$ is uniformly distributed on the circle.

Given a distribution on the line, we can wrap it around the circumference of the circle with unit radius. If $X$ is the random variable on the line, the random variable $\Theta$ of the wrapped distribution is given by
\[
\Theta=(X \mbox{ mod }2\pi ).
\] 
 A  popular example of wrapped distributions  is the wrapped Cauchy distribution with its pdf (see, for example, \citealt[p.51]{mardiajupp2000})  
\begin{equation}
f_c(\theta\vert \rho, \mu) = 
\frac{1}{2\pi}\frac{1-{\rho}^2}{1+{\rho}^2-2\rho\cos (\theta -\mu)  }, \theta\in [0,2\pi), \mu\in[0,2\pi),\rho\in [0,1),
\label{eq:wC}
\end{equation}
where $\mu$ is the mean direction parameter and $\rho$ is the concentration parameter. It is one of the  wrapped distributions whose density has a closed form and is heavy-tailed. When $\rho=0$ it also reduces to the uniform distribution. 

In what follows, Section \ref{S:construction} gives constructions of discrete circular families based on the four methods, along with some examples and results interrelating them. In particular,   we deduce families from general continuous location families, and a flexible family based on trigonometric sums. Section \ref{S:special} deduces some discrete families from established distributions such as von Mises and wrapped Cauchy.   We apply some of the models to our discrete data in Section \ref{S:examples}. In  Section \ref{S:comparisons}, we treat the problem of model misspecification for circular discrete distributions.  We conclude  the paper with a discussion in Section \ref{S:conclusion}. Some supporting material is given in the supplement.
\section{Constructions of families of circular discrete distributions}
\label{S:construction}
In this section, we elaborate on the four different methods to construct families of circular discrete distributions that were  mentioned in Section \ref{S:intro}. Although our ideas naturally extend to constructing discrete distributions on an irregular support, we will focus here on  the regular lattice support which lends itself to some mathematical simplifications. 

We denote the set of real numbers by $\mathbb{R}$, non-negative real numbers by $\mathbb{R}^+$, the set of integers by $\mathbb{Z}$, non-negative integers by $\mathbb{Z}^+$ and the cyclic group of integers modulo a given positive integer $m$ by $\mathbb{Z}_m$, i.e.
\begin{equation}
\label{eq:Zmodm}
\mathbb{Z}_m=\{0,1,\ldots, m-1\}.
\end{equation}
The regular circular lattice domain is given by the vertices of a regular polygon on the circle, denoted by $\mathcal{D}_m$, i.e.
\begin{equation}
\label{eq:circ_lattice}
\mathcal{D}_m=  \{2\pi r/m, r\in \mathbb{Z}_m \}.
\end{equation}
 We generally use $f(\cdot)$ or $g(\cdot)$ to denote a probability density function (pdf) of a continuous distribution on the line or circle,  and $p(\cdot)$ to denote a discrete probability function on $\mathbb{Z}_m$. 
\subsection{Maximum entropy discrete circular distributions}
\label{S:maxentropy}
For a probability function  $\{p(r), ~~r\in \mathbb{Z}_m\}$, with $p(r)$ denoting the probability of the point $2\pi r/m \in \mathcal{D}_m$, Shannon's entropy    is defined as
\begin{equation}
\label{eq:entropy}
~~~-\sum_{r=0}^{m-1} p(r) \log p(r).
\end{equation}
Let $t_1, t_2, \ldots, t_q$ be $q$ real valued functions defined on  $\mathbb{Z}_m$ and suppose we are interested in discrete distributions $\{p(r), ~~r\in \mathbb{Z}_m\}$ that satisfy a set of pre-selected moment conditions
\begin{eqnarray}
\label{eq:constraints_entropy}
\sum_{r=0}^{m-1}p(r) t_1(r)=a_1, ~\sum_{r=0}^{m-1}p(r)  t_2(r)=a_2, \ldots, \sum_{r=0}^{m-1}p(r)  t_q(r)=a_q,
\end{eqnarray}
with given constants $a_1, a_2, \ldots, a_q$. Then, a useful method to construct discrete distributions is to maximize the entropy among all distributions on $\mathbb{Z}_m$ that satisfy the given conditions.   As noted by \citet{kemp1997}, the philosophy behind this construction is that ``one should use all the given information and nothing else". The following theorem gives this construction, which follows from  Theorem 13.2.1  of  \citet[pp. 408-409]{kaganlinnikrao1973} on the line  and was adapted  in \citet{mardia1975} for directional distributions.
\begin{thm}[Maximum Entropy Distributions]
\label{thm:MED}
The probability function  $\{p(r), ~~r\in \mathbb{Z}_m\}$ that maximizes the entropy (\ref{eq:entropy}) subject to the constraints (\ref{eq:constraints_entropy}) is of the form
\begin{eqnarray}
 p(r)&=& \frac{ e^{\sum_{i=1}^q b_i t_i(r)}}{\sum_{k=0}^{m-1}e^{\sum_{i=1}^q b_i t_i(k)}}, ~~r\in \mathbb{Z}_m,  \label{eq:MED_gen}
 \end{eqnarray}
provided there exist constants  $b_1, b_2,\ldots, b_q$ satisfying
 \begin{eqnarray}
 \frac{  \sum_{r=0}^{m-1} t_j(r) e^{\sum_{i=1}^q b_i t_i(r)}}{\sum_{k=0}^{m-1}e^{\sum_{i=1}^q b_i t_i(k)}} =a_j, ~~~j=1,2,\ldots, q. \label{eq:MED_cond}
 \end{eqnarray}
In that case, the distribution is unique.
\end{thm}
~\\
We now give a few examples of maximum entropy discrete distributions.

\begin{eg} 
von Mises distribution: Suppose $q=2$ and $t_1(r)= \cos(2\pi r/m)$ and $t_2(r)= \sin(2\pi r/m)$, a discrete version of the von Mises distribution  is of the form
\begin{eqnarray}
 p(r) &=&  \frac{ e^{\kappa  \cos(2\pi r/m - \mu) }}{\sum_{k=0}^{m-1}e^{\kappa  \cos(2\pi r/m - \mu)}}, ~~r\in \mathbb{Z}_m,\label{eq:MED_VM}
 \end{eqnarray}
 where $\kappa = \sqrt{b_1^2+b_2^2}$ and $\tan(\mu)= b_2/b_1$. 
 ~\\
We note that this also happens to be the conditionalized discrete von Mises distribution that is discussed in more detail later.
\end{eg}
\begin{eg}
Beran distributions: A more general family than the previous example, a discrete version of the Beran family (\citealt{beran1979}), is obtained by considering constraints on the expected values of $t_{k}(r)=\left(\cos\left( 2\pi r k/m\right), \sin\left(2\pi r k/m\right) \right)$, which leads to the probability function
\begin{eqnarray}
p(r) \propto  e^{\sum_{k=1}^{q} (a_k \cos\left( 2\pi r k/m\right) + b_k \sin\left(2\pi r k/m\right))}, ~~r\in \mathbb{Z}_m. \label{eq:beran}
\end{eqnarray} 
We will denote this distribution by $\mathcal{B}_q$, where $q$ is the order of the distribution. So, $\mathcal{B}_1$ is von Mises discrete distribution as in the previous example. $\mathcal{B}_2$ is the discrete generalized von Mises distribution. 
~\\
We note that this family can be traced back to \citet{maksimov1967}, although his focus for this family is on a characterization  for the unknown centering parameter (rather than the concentration parameter), so it is of limited practical importance. 
\end{eg}
\begin{eg}
Geometric distribution: Suppose $q=1$ and $t_1(r)= r $ the maximum entropy distribution is of the form
 \begin{eqnarray}
 p(r) &= & \frac{ (1-p) p^{ r}  }{1-p^{m}}, ~~r\in \mathbb{Z}_m, \mbox{ where } p=e^{b_1}. \label{eq:MED_geom}
 \end{eqnarray}
 ~\\
Historically, \citet[p. 50]{mardia1972}  proposed the above distribution as a model for  roulette outcomes (possibly biased).
\end{eg}
The above three examples also arise out of the ``conditionalized" construction of discrete circular distributions that we  define
below in Section \ref{S:construction_marg_plug}. We note that (\ref{eq:MED_geom})  is also the ``centered wrapped geometric distribution" discussed below in the next subsection.   
\subsection{Centered wrapped discrete circular distributions}
\label{S:cent_wrap}
A natural construction  to obtain a circular discrete distribution is to start with a discrete distribution on the line and wrap it on the circle (see, for example,  \citealt[p.50]{mardia1972}). Let $Z$ be a random variable taking integer values (i.e. in $\mathbb{Z}$) with probability function $p_0(\cdot)$. For a given positive integer $m$, we define here the wrapped discrete  random variable:
$$Z_w= \left(Z\mbox{ mod } m\right)\times (2\pi/m).$$ 
We note that $Z_w\in \mathcal{D}_m$ and its probability function is given by
\begin{eqnarray}
\label{eq:pmf_wrapdisc}
p_{w0}(r)=P(Z_w=2\pi r/m)= \sum_{k=-\infty}^\infty p_0(r+km), ~~ r\in \mathbb{Z}_m.
\end{eqnarray}
It follows that the characteristic function of $Z_w$ is given by
\begin{eqnarray}
\label{eq:chfn_wrapdisc}
\psi_{p,m} = E\left( e^{i p Z_w}\right) = \phi ( 2\pi p/m),
\end{eqnarray}
where $\phi(\cdot)$ is the characteristic function of $Z$. In general, these distributions do not have a  mean direction or centering parameter, and therefore we  construct the ``centered wrapped" probability function with a centering  parameter $t$ as follows,
\begin{eqnarray}
p(r)= \begin{cases} p_{w0}(r-t+m), ~~ r<t  \\ p_{w0}(r-t),  ~~~~~~~~ r\ge t, \end{cases} ~~ r, t\in \mathbb{Z}_m.
\end{eqnarray}
Choosing the domain of $t$ as $\mathbb{Z}_m$ ensures probabilities are well defined without changing the domain of the distribution.  
\begin{eg}
Centered wrapped Poisson distribution: For the Poisson distribution with mean $\lambda$, the wrapped Poisson distribution has the probability function
\begin{eqnarray}
\label{eq:wrap_poi_pw0}
p_{w0}(r) = e^{-\lambda}\sum_{k=0}^{\infty} \frac{\lambda^{r+km}}{(r+km)!}, ~r\in \mathbb{Z}_m.
\end{eqnarray}
The centered wrapped probability function with centering parameter $t$ is then given by
 \begin{eqnarray}
\label{eq:wrap_poi_pw}
p(r) = \begin{cases} e^{-\lambda}\sum_{k=0}^{\infty} \frac{\lambda^{r-t+m+km}}{(r-t+m+km)!}, ~~r< t \\ e^{-\lambda}\sum_{k=0}^{\infty} \frac{\lambda^{r-t+km}}{(r-t+km)!}~~~, ~~r\geq  t\end{cases} ~r,t \in \mathbb{Z}_m.
\end{eqnarray}
The above expression is a special case of the distributions considered by \citet{Mastrantonioetal2019}. 
\end{eg}
For practical applications with continuous circular data,  it is well known that the selected  probability density is continuous at $2\pi$, i.e. the pdf value at $0$ is same as its limiting value at $2\pi$. Similarly, along the same lines, a 
 desirable property for a circular discrete probability function $p(\cdot)$ on  $\mathbb{Z}_m$ is to have $p(0)=p(m)$. The maximum entropy  and the centered wrapping methods do not necessarily ensure this property as is apparent from Examples 3 and 4. However, this property is ensured if we construct discrete distributions by applying the marginalized and conditionalized  methods on continuous circular distributions, which we will discuss next.
\subsection{Marginalized and conditionalized discrete  distributions}
\label{S:construction_marg_plug}
In this section,  we focus on univariate circular constructions based on marginalized and conditionalized methods, whose brief descriptions were given in Section \ref{S:intro}. 

There has been literature on the marginalized and conditionalized discretization on the line. For example \citet{kemp1997} and \citet{szablowski2001}  discuss the conditionalized discrete normal, \citet{inusaha_kozubowski2006} discuss the conditionalized discrete Laplace, and \citet{papadatos2018} derives the characteristic function of the conditionalized discrete Cauchy. The conditionalized approach can also be described as a ``plug-in" approach, whereas the marginalized approach is in fact equivalent to the well known probit construction, usually used for univariate and multivariate normal  distributions, see for example \citet[p.20]{joeharry2014}.  \cite{Alzaatrehetal2012}  and \cite{Chakraborty2015} make references to both of these methods, while developing other methods of constructions. However, there have not been any insights relating these constructions.

Marginalized and conditionalized constructions of circular families of discrete distributions are very recent as proposed in  \citet{mardiasriram2020ar}. Besides, there has not been published in-depth analysis of truly discrete circular data. However, particular cases of the conditionalized approach, including von Mises and wrapped Cauchy distributions, have appeared, not only in \citet{mardiasriram2020ar} but also in \cite{girijaetal2019} and \citet{imotoetal2020}.  We give a unified treatment of the different methods as a strategy to construct rich classes of discrete distributions on the circle.  We derive new  results (e.g. Theorems \ref{thm:MED} to \ref{thm:duality}) that offer insights on the inter-relationships between the constructions. 

We now define the marginalized and conditionalized discrete families for  the circular case. There has been some very recent work on these approaches although not in a comprehensive and unified way, as we describe below. Let $\Theta$ be a random variable with pdf $f(\theta), ~\theta\in[0, 2\pi)$. 
\begin{defin}
\label{def:MD_circ}
 The probability function of the marginalized discrete (MD) distribution on the circle is given by
\begin{equation}
\label{eq:pmf_by_modulo}
p(r)=  \int_{\frac{2\pi r}{m}}^{\frac{2\pi (r+1)}{m}} f(\theta) d\theta= F\left(\frac{2\pi (r+1)}{m}\right) - F\left(\frac{2\pi r}{m}\right), ~r\in\mathbb{Z}_m,
\end{equation}
where $F(\cdot)$ is the cumulative distribution function of the pdf $f(\cdot)$. 
\end{defin}
We note that this is also the probability function of the discrete random variable  $\lfloor \frac {m \Theta} {2 \pi}\rfloor$, where $\lfloor \cdot \rfloor$ denotes the largest integer less than or equal to the given number.
\begin{defin}
\label{def:CD_circ}
The probability function of the conditionalized discrete (CD) distribution on the circle is given by
\begin{equation}
\mbox{p} \left( r \right) = \frac{f (2\pi r/m)}{\sum_{k=0}^{m-1} f(2\pi k/m)},  r\in\mathbb{Z}_m.
\label{eq:plugin_method} 
\end{equation}
\end{defin}
For simplicity, we will denote both probability functions (\ref{eq:pmf_by_modulo}) and (\ref{eq:plugin_method}) by the same notation, but the choice will be obvious from the context. One question of interest is ``can the marginalized and conditionalized methods   lead to the same discrete distribution on the circle?"  We show that under certain conditions, the two approaches will lead to diffent discrete distributions on the circle except for the trivial uniform case. This has implication when we come to selecting between the two approaches in practice and we give insights based on comparison of the two approaches for some particular cases in Section \ref{S:comparisons}.  Theorem \ref{thm:char_on_circle} below gives one characterization. This theorem is inspired by a similar question for the linear case related to the exponential distribution (see supplement).
\begin{thm}
\label{thm:char_on_circle}
Suppose $f(\cdot)$ is a strictly positive and continuous circular pdf on $[0, 2\pi] $ with $f(\theta)=f(\theta+2\pi)$. Then,  the two discretization approaches (i.e. MD and CD) lead to the same discrete distribution, with \begin{eqnarray}
\label{eq:margeqcond_circ}
\frac{f(a+k\delta)}{\sum_{r=0}^{m-1} f(a+r\delta)}&&= \frac{\int_{a+k\delta}^{a+(k+1)\delta} f(\theta) d\theta}{\int_{a}^{a+m\delta} f(\theta) d\theta},  \forall ,  a, \delta\in[0, 2\pi), ~~\forall~k\in \mathbb{Z}_m \mbox{ with }  m\delta\leq 2\pi,
\end{eqnarray}
 iff $f$ is the uniform density.
  \begin{proof}
Considering (\ref{eq:margeqcond_circ})  for $k=1$ and $k=0$, and taking their ratio, we get 
\begin{eqnarray}
f(a+\delta)\int_{a}^{a+\delta}f(\theta) d\theta =f(a)\int_{a+\delta}^{a+2\delta}f(\theta) d\theta. \label{eq:char_circ1}
\end{eqnarray}
Integrating both the left hand side (lhs) and right hand side (rhs) of (\ref{eq:char_circ1}) with respect to $\delta\in[0, 2\pi)$, we get
\begin{eqnarray}
\int_{0}^{2\pi}f(a+\delta)\int_{a}^{a+\delta}f(\theta) d\theta d\delta 
=f(a) \int_{0}^{2\pi}\int_{a+\delta}^{a+2\delta}f(\theta) d\theta d\delta. \label{eq:char_circ2}
\end{eqnarray}
Using continuity of $f(\cdot)$ including $f(\theta+2\pi)= f(\theta)$, we have $\int_{a}^{a+2\pi}f(\theta) d\theta=1$ for any $a$, and the lhs of (\ref{eq:char_circ2}) can be simplified as
\begin{eqnarray}
\int_{0}^{2\pi}f(a+\delta)\int_{a}^{a+\delta}f(\theta) d\theta d\delta= \frac{1}{2} \left( \int_{a}^{a+\delta}f(\theta) d\theta\right)^2 \vert_{\delta=0}^{2\pi}= 1/2. \label{eq:char_circ3}
\end{eqnarray}
Now, the rhs of (\ref{eq:char_circ2}) can be shown to be
\begin{eqnarray}
f(a) \int_{0}^{2\pi}\int_{a+\delta}^{a+2\delta}f(\theta) d\theta d\delta = f(a) (A- B) = f(a) \pi, \label{eq:char_circ4}
\end{eqnarray}
Equation (\ref{eq:char_circ4}) follows because $A$ and $B$ can be simplified as below. 
\begin{eqnarray*}
A&=& \int_{0}^{2\pi}\int_{a}^{a+2\delta}f(\theta) d\theta d\delta = \delta \int_{a}^{a+2\delta}f(\theta) d\theta  \vert_{0}^{2\pi} - \int_0^{2\pi}2\delta f(a+2\delta) d\delta\\
&=& 2\pi \int_{a}^{a+4\pi}f(\theta) d\theta -\frac{1}{2}\int_{0}^{4\pi}\delta^\prime f(a+\delta^\prime) d\delta^\prime=4\pi - \left(\pi + \int_{0}^{2\pi}\delta^\prime f(a+\delta^\prime) d\delta^\prime\right),
\end{eqnarray*}
and 
\begin{eqnarray*}
B&=& \int_{0}^{2\pi}\int_{a}^{a+\delta}f(\theta) d\theta d\delta=\delta \int_{a}^{a+\delta}f(\theta) \vert_0^{2\pi} - \int_0^{2\pi} \delta f(a+\delta)d\delta = 2\pi - \int_0^{2\pi} \delta f(a+\delta)d\delta.
\end{eqnarray*}
Equating the lhs (\ref{eq:char_circ3}) and rhs (\ref{eq:char_circ4}) , we get $f(a) = \frac{1}{2\pi}.$
Since $a$ is arbitrary, this  means that $f(\cdot)$ must be the uniform pdf on the circle.
     \end{proof}
\end{thm}
 It is to be emphasized that the assumption on continuity of the pdf $f(\cdot)$ is crucial in the above theorem.  For example, the marginalized and conditionalized methods applied to the wrapped exponential distribution on the circle lead to the same discrete distribution, i.e. geometric distribution. However, the wrapped exponential pdf is not continuous at $\theta=2\pi$.  

A related question is whether the marginalized and conditionalized methods  can belong to the same family of distributions. Indeed, it is easy to see that this property will hold for a generalized  Cardioid-type family of distributions as given by the following theorem. 
\begin{thm}~\\
\label{eq:pD_MD_circ_faminv}
Consider the pdf of the parent family  defined by
\begin{eqnarray}
\label{eq:pdf_cardioid_mix}
f(\theta) =  \frac{1}{2\pi}\sum_{k=1}^\infty \eta_k\left(1+2\rho_k \cos(\theta - \mu_k)  \right), ~~\theta\in[0, 2\pi),
\end{eqnarray}
where $\sum_{k=1}^\infty \eta_k=1$, $\forall ~k$, $\eta_k\geq 0,~\mu_k\in ~[0, 2\pi), ~|\rho_k|<1/2$. Then the  marginalized discrete distribution   is also a member of   the family of conditionalized discrete distributions.
\end{thm}

An interesting connection between the constructions on the circle and line, is given by the following theorem.
\begin{thm}[Duality]
\label{thm:duality} ~\\ Consider the following dual approaches to constructing discrete circular distributions supported on $ \mathbb{Z}_m$, starting with a real valued random variable $X$, with a pdf $f(\cdot)$ on $\mathbb{R}$,  via either the marginalized  or the conditionalized methods of discretization.
\begin{itemize}
\item[(a)] {\bf Scale, discretize and wrap:} Start with the pdf of the scaled random variable $\widetilde{X}= mX/(2\pi)$, obtain the marginalized [or  conditionalized] discrete probability function on the line and denote the corresponding random variable by $\widetilde{X}_d$. Further, wrap $\tilde{X}_d$, i.e. let  $\widetilde{X}_{dw}=(\widetilde{X}_d\mbox{ mod }m)$.
\item[(b)] {\bf Wrap, scale  and discretize:} Let $X_w=(X \mbox{ mod }2\pi )$ (i.e.  $X$ wrapped on the circle). Now, start with the pdf of the scaled random variable $\widetilde{X}_w= mX_w/(2\pi)$, obtain the marginalized [or conditionalized] discrete probability function on the circle and denote the corresponding random variable by $\widetilde{X}_{wd}$.
\end{itemize}
Then,  $\widetilde{X}_{dw}$ and $\widetilde{X}_{wd}$ have the same distribution.
\begin{proof}
For a given random variable $Z$ with continuous support, recall that discretization by marginalized method means taking the random variable $\lfloor Z \rfloor$, and discretization by conditionalized method means taking  the random variable $Z_d$ with its probability function defined by $P(Z_d=r)=\frac{f(r)}{\sum_{k}f(k)}$. 
First, we will prove the equivalence for the conditionalized discretization method. The probability function of the discrete distribution resulting from process (a) is given by
\begin{eqnarray}
P(\widetilde{X}_{dw}= r)&=& \sum_{k\in \mathbb{Z}}P(\widetilde{X}_d= r+km) =
\frac{\sum_{k\in \mathbb{Z}} f(2\pi(r+km)/m)}{\sum_{k\in \mathbb{Z}}f(2\pi k/m)} .\label{eq:plug_a}
\end{eqnarray}
For the process in (b), let us denote the pdf of $X_w$ by $f_{x_w}$. Then, the probability function of the discrete distribution resulting from process (b) will be  
\begin{eqnarray}
P(\widetilde{X}_{wd}= r)= \frac{f_{x_w}(2\pi r/m)}{\sum_{r=0}^{m-1}f_{x_w}(2\pi r/m)} 
= \frac{\sum_{k\in \mathbb{Z}} f(2\pi r/m + 2\pi k)}{\sum_{k\in \mathbb{Z}} f(2\pi k/m )}. \label{eq:plug_b}
\end{eqnarray}
Since (\ref{eq:plug_a}) and (\ref{eq:plug_b}) are the same, (a) and (b) yield the same discrete circular distribution. Now, to prove the equivalence under the marginalized method of discretization, we observe that process (a) leads to the probability function given by
\begin{eqnarray}
P(\widetilde{X}_{dw}=r) &=&
\sum_{k\in \mathbb{Z}} P(\widetilde{X}_d = r +km )
= \sum_{k \in \mathbb{Z}} \left( F\left( \frac{ 2\pi (r+1)}{m} + 2\pi k \right) - F\left( \frac{2\pi r}{m} + 2\pi k\right) \right),\nonumber \label{eq:marg_a}
\end{eqnarray}
but process (b) also leads to the same probability function because
\begin{eqnarray}
P(\widetilde{X}_{wd}=r) &=& P\left(X_w \in \left[\frac{2\pi r}{m}, \frac{2\pi (r+1)}{m}\right]\right)=\sum_{k \in \mathbb{Z}} \left( F\left( \frac{ 2\pi (r+1)}{m} + 2\pi k \right) - F\left( \frac{2\pi r}{m} + 2\pi k\right) \right).\nonumber \label{eq:marg_b}
\end{eqnarray}
\end{proof}
\end{thm}
\subsubsection{General circular discrete location family}
Consider a general circular location family (see, for example,  \citealt{mardia1975b}) with probability density function given by 
\begin{equation}
\label{eq:loc_fam}
f(\theta\vert\tau, \mu)= g_\tau(\theta - \mu), ~~\theta, \mu\in[0, 2\pi), ~\tau\geq 0 ,
\end{equation}
which we assume to be unimodal with mode at $\mu$. 
For simplicity, we assume $g_\tau(\theta)=g_\tau(2\pi- \theta)$,  $g_\tau(\theta)>0$ for all $\theta\in [0, 2\pi)$ and also that $g_\tau(2\pi)= g_\tau(0)$. Note that the normalizing constant will depend only on $\tau$ and not on $\mu$. Here,
$\tau\geq 0$ is another parameter in addition to $\mu$, such that  $\tau=0$ corresponds to the case of uniform distribution and the dispersion around the mode decreases as $\tau$ increases. For example,  $\tau=\kappa$ for the von Mises distribution (\ref{eq:von_Mises1}), and $\tau=\rho$ for wrapped Cauchy distribution (\ref{eq:wC}).    

The probability function for the marginalized discrete distribution based on the circular location family (\ref{eq:loc_fam}) is given by
\begin{eqnarray}
\label{eq:loc_fam_marg}
p(r\vert m, \tau, \mu)= \int_{2\pi r/m}^{2\pi (r+1)/m} g_\tau(\theta- \mu) d\theta, ~~ r\in \mathbb{Z}_m, ~\mu\in[0, 2\pi),
\end{eqnarray}
and we will call this distribution the ``marginalized discrete circular location family". 

Similarly, the probability function of the conditionalized discrete distribution based on the circular location family (\ref{eq:loc_fam}) is given by
\begin{equation}
\label{eq:PDgen}
p(r\vert m, \tau ,\mu)= \frac{g_\tau(2\pi r/m-\mu )}{\sum_{r=0}^{m-1} g_\tau(2\pi r/m -\mu)}, ~r\in \mathbb{Z}_m, \mu\in[0,2\pi),
\end{equation}
and we will call this distribution the ``conditionalized discrete circular location family". 

The characteristic function for the probability function (\ref{eq:loc_fam_marg}) or (\ref{eq:PDgen}) is given by
\begin{equation}
\label{eq:genchfn}
\psi_{p,m}= \sum_{r=0}^{m-1}p(r\vert m, \tau ,\mu) e^{ip\frac{2\pi r}{m}}.
\end{equation}
Since $g_\tau(\theta)$  can be expressed in terms of its characteristic function  $\phi_p$ (see \citealt[p.27]{mardiajupp2000}) as 
$$ g_{\tau}(\theta)= \frac{1}{2\pi}\sum_{q=-\infty}^{\infty} \phi_q e^{-iq\theta},$$ it can shown that the characteristic function of the marginalized discrete location family is
\begin{eqnarray}
\psi_{p,m} &=& \begin{cases} 1, ~~~ p=0 \\ \frac{me^{ip(\mu-\frac{\pi}{m})}\sin(\pi p/m)}{\pi} \sum_{l=-\infty}^{\infty} \frac{\phi_{lm+p}}{lm+p}e^{i lm \mu}, p\in \mathbb{Z}_m\backslash \{0\},\end{cases} \label{eq:gen_chfn_MD}
\end{eqnarray}
and the characteristic function of the conditionalized discrete location family is
\begin{eqnarray}
\psi_{p,m} &=&  e^{ip\mu}\frac{\sum_{l=-\infty}^\infty \phi_{lm+p}e^{i lm\mu}}{\sum_{l=-\infty}^\infty \phi_{lm}e^{ilm\mu}},~p\in\mathbb{Z}_m. \label{eq:gen_chfn_CD}
\end{eqnarray}

\subsubsection{Circular discrete family based on trigonometric  sums}
\label{S:MDCDFS}
In this section, we derive the marginalized and conditionalized discrete distributions starting from the flexible continuous distribution based on trigonometric sums introduced by \citet{fernandez2004}.  For a set of complex numbers  ${\bf c}= \{c_0, c_1,\ldots, c_J\}$  such that
\begin{equation}
\sum_{k=0}^{J} |c_k|^2=\frac{1}{2\pi}, \label{eq:FD_ck}
\end{equation}
 the pdf defined by \citet{fernandez2004}  is 
\begin{eqnarray}
f(\theta)&= &  \frac{1}{2\pi}+ \frac{1}{\pi}\sum_{k=1}^{J} \left\{a_k \cos(k\theta) + b_k \sin(k\theta)\right\}, ~~\theta\in[0, 2\pi), \label{eq:FD_f}
\end{eqnarray}
where $(a_k,b_k)$ are such that $a_k-i b_k = 2\sum_{\nu=0}^{J-k} c_{\nu+k} \bar{c}_\nu$.    The specific choice of  $(a_k, b_k)$  is a necessary and sufficient condition to ensure positivity of  the function $f$. The pdf (\ref{eq:FD_f}) can also be written as
\begin{eqnarray}
f(\theta)&= &  \frac{1}{2\pi}+ \frac{1}{\pi}\sum_{k=1}^{J} \rho_k \cos(k\theta-\phi_k), ~~\theta\in[0, 2\pi), \label{eq:FD_f2}
\end{eqnarray}
where 
\begin{equation}
\label{eq:FDrhok}
\rho_k=\sqrt{a_k^2+b_k^2} \mbox{ and } \phi_k=\arctan(b_k/a_k).
\end{equation}

We will refer to the distributions obtained by applying Definitions  \ref{def:MD_circ} and  \ref{def:CD_circ} on the pdf (\ref{eq:FD_f}) as the ``marginalized discrete trigonometric sum"  distribution (denoted $MDTS(m, \boldsymbol{c})$) and the ``conditionalized discrete trigonometric sum" distribution (denoted $CDTS(m, \boldsymbol{c})$), respectively.  It is easy to see from (\ref{eq:FD_f}) and (\ref{eq:FD_f2}) that the probability function of $MDTS(m, \boldsymbol{c})$  is given by
\begin{eqnarray}
&&p(r|m,\boldsymbol{c})=  \frac{1}{m}+\frac{2}{\pi} \sum_{k=1}^J \frac{\sin(\pi k/m)\rho_k}{k} \cos\left( \frac{2\pi k(r+1/2)}{m} - \phi_k\right) , ~r \in \mathbb{Z}_m,\label{eq:MDTS}
\end{eqnarray}
and the probability function of $CDTS(m, \boldsymbol{c})$ is given by
\begin{eqnarray}
p(r|m,\boldsymbol{c})= \frac{1+2\sum_{k=1}^J \rho_k \cos\left(\frac{2\pi k r}{m}-\phi_k\right)}{m\left(1+2\underset{k=0 (\mbox{mod m})}{\sum_{k=1}^J }a_k\right)}, ~r \in \mathbb{Z}_m,  \label{eq:CDTS}
\end{eqnarray}
where $\rho_k$ and $\phi_k$ are as in equation (\ref{eq:FDrhok}). Note that for $J=1$, equations (\ref{eq:MDTS}) and (\ref{eq:CDTS}) give the marginalized and conditionalized discrete cardioid distributions, respectively. We see that the two distributional forms are identical although with different parametrization, which is consistent with Theorem \ref{eq:pD_MD_circ_faminv}.

 As in the continuous case, the above constructed discrete families give flexibility in modeling multimodality and skewness in circular discrete data. We note that the probability function for CDTS is derived in \citet{imotoetal2020}. In this paper, we will give a particular application in Section \ref{S:roulettefixedframe}.

\section{Key special discrete distributions and their properties }
\label{S:special}
We now give  the marginalized and conditionalized methods for   the von Mises (\ref{eq:von_Mises1}) and the wrapped Cauchy (\ref{eq:wC}) as the parent distributions followed by  some basic properties including characteristic function, estimation and hypothesis testing. We begin with the definitions of these distributions.\\
   
\begin{defin}
Using (\ref{eq:pmf_by_modulo}), the  probability function for the marginalized discrete von Mises  distribution (MDVM) with mean parameter $\mu$ and concentration parameter $\kappa$, denoted by $MDVM(m, \kappa, \mu)$, is given by
 \begin{equation}
 \label{eq:PDVM_0md}
 p(r\vert m, \kappa, \mu) = \frac{1}{2\pi I_0(\kappa)}\int_{2\pi r/m}^{2\pi(r+1)/m}e^{\kappa \cos\left(\theta-\mu\right)}d\theta, r\in\mathbb{Z}_m,~~\mu\in[0, 2\pi),
 \end{equation}
\end{defin}

\begin{defin}
 Using (\ref{eq:plugin_method}), the  probability function for the conditionalized discrete von Mises (CDVM) distribution with mean parameter $\mu$ and concentration parameter $\kappa$, denoted by $CDVM(m, \kappa, \mu)$, is given   by
 \begin{equation}
 \label{eq:PDVM_0cd}
 p(r\vert m, \kappa, \mu) = \frac{1}{L_0(\kappa,\mu)}e^{\kappa \cos\left(2\pi r/m-\mu\right)}, r\in\mathbb{Z}_m,~~\mu\in[0, 2\pi),
 \end{equation}
where the normalizing constant  is the reciprocal of the function
 \begin{equation}
 \label{eq:C0_CDVM}
 L_0(\kappa, \mu)= \sum_{r=0}^{m-1} e^{\kappa \cos\left( 2\pi r/m-\mu\right)}.
 \end{equation}
 \end{defin}
Similarly, we have  the following definitions for the wrapped Cauchy case.

\begin{figure}[h]
\center
\begin{subfigure}{.48\textwidth}
\caption{~}
\includegraphics[width=5cm,angle=270]{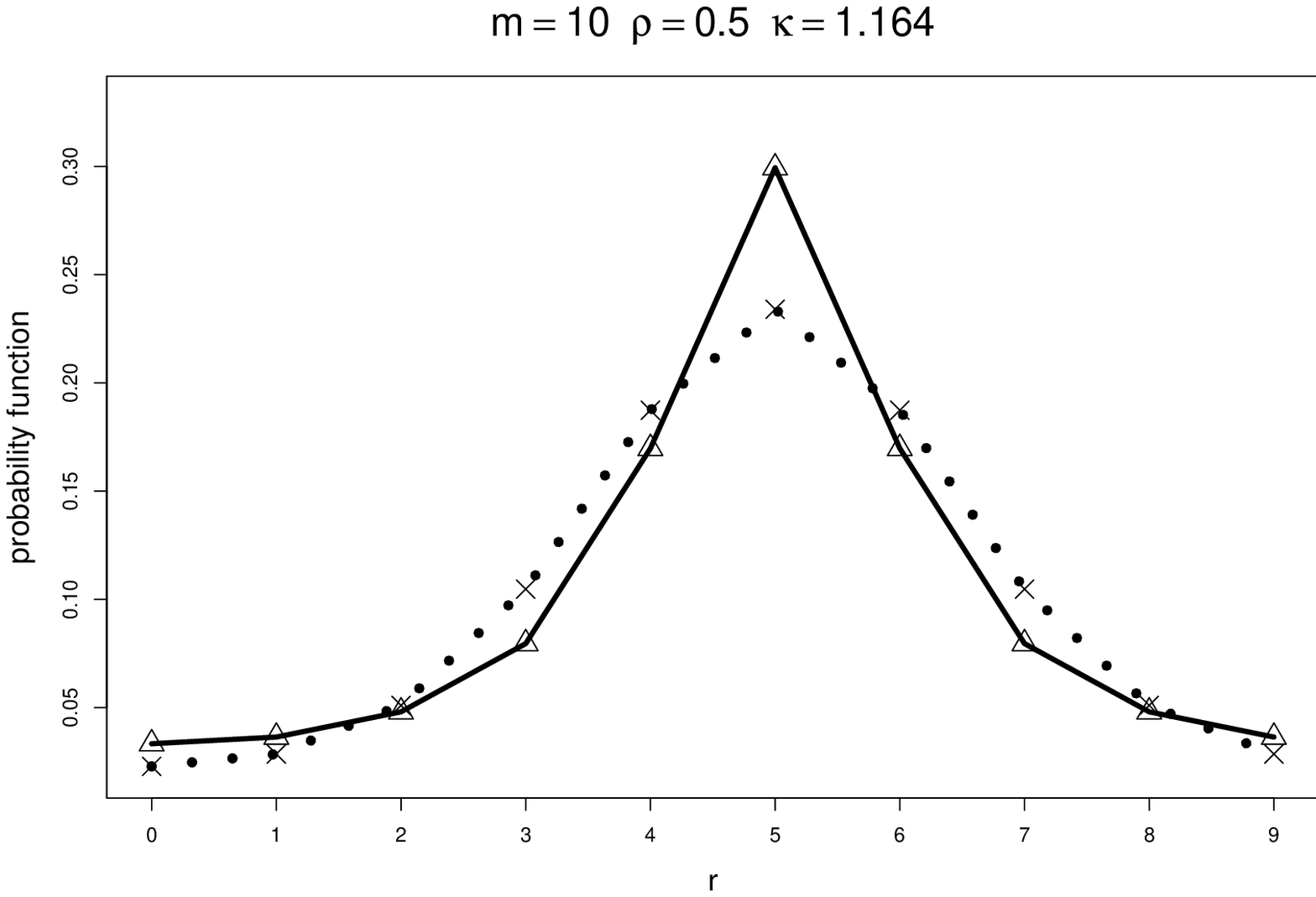}
\end{subfigure}
\begin{subfigure}{.48\textwidth}
\caption{~}
\includegraphics[width=5cm,angle=270]{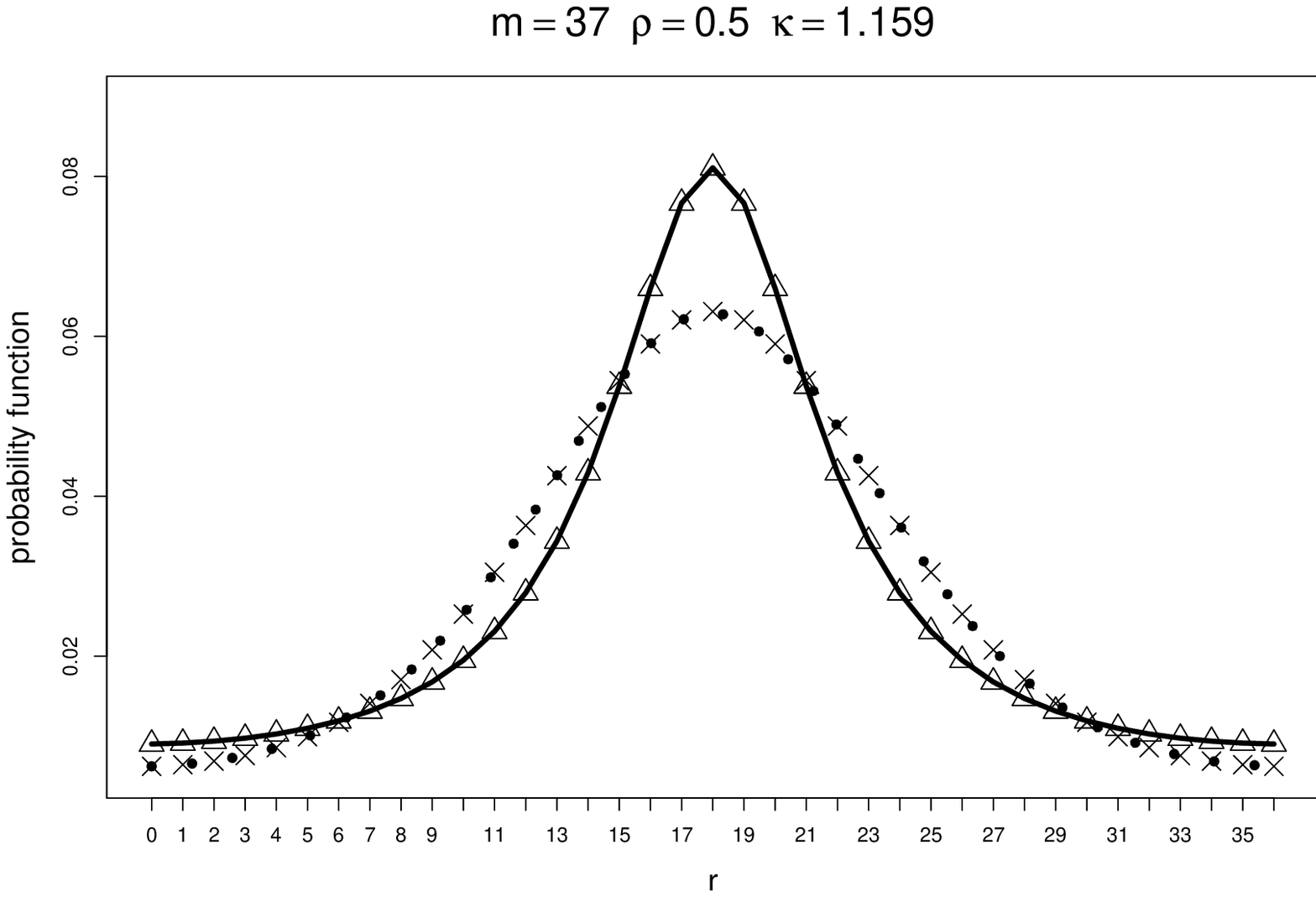}
\end{subfigure}
\caption{\label{F:pdwc_pdvm_pmf00} Probability  functions of $CDWC(m, \rho, \mu)$(triangles joined by solid line) and $CDVM(m, \kappa, \mu)$ (cross joined by dotted lines) plotted for (i) $m=10$ and (ii) $m=37$ with  $\mu=2\pi 5/m$ and $\mu=2\pi 16/m$ respectively, for $\rho=0.5$ and its mapped $\kappa$ value by matching the first trigonometric moment.}
\end{figure}

\begin{defin}
The probability function of the marginalized discrete wrapped Cauchy (MDWC) distribution with mean parameter $\mu$ and concentration parameter $\rho$, denoted by $MDWC(m, \rho, \mu)$, is given by
\begin{eqnarray}
p(r\vert m, \rho, \mu)&=& \frac {1}{2\pi} 
\cos ^{-1} \left\{ \frac {(1 + \rho^2) \cos ( \frac{2\pi (r+1)}{m} - \mu ) - 2\rho}
		{1 + \rho^2 - 2\rho \cos (\frac{2\pi (r+1)}{m} - \mu )} \right\} \nonumber\\
		&~& - \frac {1}{2\pi} 
\cos ^{-1} \left\{ \frac {(1 + \rho^2) \cos ( \frac{2\pi r}{m} - \mu ) - 2\rho}
		{1 + \rho^2 - 2\rho \cos (\frac{2\pi r}{m} - \mu )} \right\} \nonumber\\
		&& ~~~~~~r\in\mathbb{Z}_m, ~~\mu\in [0, 2\pi), ~~\rho\in[0,1).
\end{eqnarray}
\end{defin}
Alternatively, for computational purposes, the above expression can be written as 
\begin{eqnarray}
p(r\vert m, \rho, \mu)&=& \frac{1}{\pi} \arctan\left(\frac{ \frac{1+\rho}{1-\rho}\left\{\tan(\pi(r+1)/m- \mu/2)- \tan(\pi r/m-\mu/2)\right\}}{1+ \left(\frac{1+\rho}{1-\rho}\right)^2\tan(\pi(r+1)/m-\mu/2)\tan(\pi r/m -\mu/2) }\right),~  \nonumber \\
&& ~~r\in\mathbb{Z}_m, ~~\mu\in [0, 2\pi), ~~\rho\in[0,1). \label{eq:MDWC_0}
\end{eqnarray}
Further, we have  
\begin{defin}
 The probability function of the conditionalized discrete wrapped Cauchy  distribution with mean parameter $\mu$ and concentration parameter $\rho$, denoted by $CDWC(m, \rho, \mu)$, is given by
\begin{eqnarray}
p(r\vert m, \rho, \mu)&= &\frac{1}{D_0(\rho, \mu)}\frac{1}{1-2\rho \cos\left(\frac{2\pi r}{m}-\mu\right) +\rho^2},\nonumber \\
&&~~ ~r\in \mathbb{Z}_m,~ \mu\in[0, 2\pi),~ \rho\in[0,1).\label{eq:PDWC_mu}
\end{eqnarray}
where the normalizing constant  is the reciprocal of the function
 \begin{equation}
 \label{eq:C0_CDWC}
 D_0(\rho, \mu)= \frac{m(1-\rho^{2m})}{(1-2\rho^m \cos(m\mu)+\rho^{2m})(1-\rho^2)}.
 \end{equation}
\end{defin}
The normalizing constant (\ref{eq:C0_CDWC})  is derived in the supplement. 
 For simplicity of notation, while writing the marginalized or conditionalized discrete probability functions, we  will omit the subscripts such as $v$ or $c$~corresponding to von Mises or Cauchy, but they will be clear from the context and by the explicit mention of $\kappa$ versus $\rho$ as the concentration parameters.  So, we will always denote the discrete probability functions by  $p(r \vert m,\kappa, \mu)$ for von Mises or $p(r \vert m, \rho, \mu)$ for wrapped Cauchy.
 We now make some additional observations specific to  CDVM and CDWC distributions.  \\ 
~\\
{\bf Probability functions.} Figure \ref{F:pdwc_pdvm_pmf00} plots the probability  functions of $CDWC(m, \rho, \mu)$ and $CDVM(m, \kappa, \mu)$, for (i) $m=10$ and (ii) $m=37$ with  $\mu=2\pi 5/m$ and $\mu=2\pi 16/m$ respectively, for $\rho=0.5$ and its mapped $\kappa$ value.  In order to compare the probability functions of $CDWC(m, \rho, \mu)$ and $CDVM(m, \kappa, \mu)$, we need to first map the parameters $\rho$ to $\kappa$. We do so by matching their first trigonometric moments given by equations (\ref{eq:PDVM8}) and (\ref{eq:rho_w}) below, i.e. $ B(\kappa)=\rho_w$. We note that the CDWC is more spiked and heavy tailed compared to CDVM. \\
~\\

~\\

{\bf Characteristic functions.} We make some observations based on the characteristic functions of $CDVM$ and $CDWC$ distributions.\\
{\bf (a) CDVM.} For the $CDVM$ distribution (with $\mu=0$), we can also obtain an alternative simplified form for the characteristic function.  Let us write
\begin{equation}
\label{eq:pdvm1}
L_{p}(\kappa)= \sum_{r=0}^{m-1}\cos\left(p\frac{2\pi r}{m}\right)e^{\kappa \cos\left(\frac{2\pi r}{m}\right)}.
\end{equation}
So, for $\Theta \sim CDVM(m, \kappa, \mu=0)$, we have
\begin{eqnarray}
\psi_{p,m}=E\left[e^{i p \Theta} \right]&=&  B_p(\kappa) , \mbox{where }B_{p}(\kappa)= L_{p}(\kappa)/L_{0}(\kappa).\label{eq:chfn_PDVM}
\end{eqnarray}
It then follows by writing $B_1(\kappa)=B(\kappa)$, that 
\begin{equation}
B(\kappa)= E\left(\cos \frac{2\pi r}{m}\right) \mbox{ and } B^\prime(\kappa)= Var\left(\cos \frac{2\pi r}{m}\right) \label{eq:PDVM8}.
\end{equation}
$L_p(\kappa)$ is the discrete analogue of the Bessel function $I_p(\kappa)$. It is to be noted that not all of the standard identities of $I_p(\kappa)$ (see \citet[Appendix 1]{mardiajupp2000} necessarily hold for its discrete analogue.\\
{\bf (b) CDWC.}  For $\Theta \sim CDWC(m, \rho, \mu=0)$, the $p$th trigonometric moments ($p\in \mathbb{Z}_m$) are given by
\begin{equation}
\label{eq:beta_p}
 \alpha_{p,m}=E\left[\sin \left(p\Theta\right)\right]=0, ~~\beta_{p,m}=E\left[\cos \left(p\Theta\right) \right]= \frac{\rho^p(1+\rho^{m-2p})}{1+\rho^m}.
\end{equation}
For $p=1$, this leads to the mean resultant length
\begin{equation}
\label{eq:rho_w}
\rho_w = \frac{\rho (1+\rho^{m-2})}{1+\rho^m}.
\end{equation}
 In general, $0\leq \rho_w \leq 1$ and as $m\rightarrow \infty$, $\rho_w \rightarrow \rho$. By property of characteristic functions, in general $0\leq \beta_{p,m}\leq 1$, and as $m\rightarrow \infty$, $\beta_{p,m} \rightarrow \rho^p$, as can be seen from Equation (\ref{eq:beta_p}), which is known to be the characteristic function of the wrapped Cauchy distribution  as expected. Further, this convergence happens at an exponential rate. To see this, note  that
\begin{eqnarray*}
\psi_{p,m}(\rho) - \rho^p &=& \frac{\rho^p \left(\rho^{m-2p}- \rho^m\right)}{1+\rho^m}= \frac{ \rho^{m-p} \left(1-\rho^{2p}\right)}{1+\rho^m}.
\end{eqnarray*}
For any fixed $p$, it follows that $|\psi_{p,m}(\rho) - \rho^p| \leq \rho^{m-p}$, and hence $|\psi_{p,m}(\rho) - \rho^p| = \mathcal{O}(\rho^{m-p})$.  In particular, since $\psi_{1,m}= \rho_w$, we have  $|\rho_w - \rho| = \mathcal{O}(\rho^{m-1})$.  

{\bf Estimation.} The maximum likelihood estimates (mle) of $(\mu, \kappa)$ for $CDVM$ and $(\mu, \rho)$ for $CDWC$ can be obtained by iteratively solving the maximum likelihood equations for the two parameters, the details of which are given in the supplement. Further, we give  asymptotically equivalent  estimates to mle which are simpler to compute:\\
{\bf (a).} Moment estimates: Characteristic functions for the marginalized discrete and conditionalized discrete location families with cardioid, von Mises and wrapped Cauchy as parent distributions, based on the general formulas (\ref{eq:gen_chfn_MD}) and (\ref{eq:gen_chfn_CD}),  are given in the supplement. 
These can be used to estimate parameters based on matching of trigonometric moments from the data. For $\Theta \sim CDWC(m, \rho, \mu)$, the trigonometric moments have a closed form,  given by
\begin{eqnarray*}
E[\cos(\Theta)]&=& A\cos(\mu) + B\cos((m-1)\mu),\\ 
E[\sin(\Theta)]&=& A\sin(\mu) - B\sin((m-1)\mu),
\end{eqnarray*}
where $$A=\frac{\rho(1-\rho^{2m-2})}{(1-\rho^{2m})}, B=\frac{\rho^{m-1}(1-\rho^2)}{(1-\rho^{2m})}.$$
For the constrained case when $\mu = 2\pi t/m, t\in \mathbb{Z}_m$, the above equations simplify to
$$E[\cos(\Theta)]=\frac{\rho (1+\rho^{m-2})}{1+\rho^m}\cos(\mu) \mbox{ and }E[\sin(\Theta)]=\frac{\rho (1+\rho^{m-2})}{1+\rho^m}\sin(\mu). $$
{\bf (b).} Hybrid estimates: When $n$ is large, plug-in the sample mean direction ($\bar{\theta}$) for $\mu$, and obtain mle for only $\rho$ or $\kappa$.\\
{\bf (c).} Constrained estimates: When $m$ is large, constrain $\mu$ to $\{2\pi t/m, t\in \mathbb{Z}_m\}$. In this case, the normalizing constants (\ref{eq:C0_CDVM}) and (\ref{eq:C0_CDWC}) become free of $\mu$,  leading to some simplifications.  For $CDVM$,  this approach leads to following explicit equations for  mle:
$$B(\hat{\kappa}) = \bar{R} \cos\left(\bar{\theta} - \frac{2\pi \hat{t}}{m}\right) \mbox{ and }~~\hat{t}= \left[\frac{m\bar{\theta}}{2\pi}\right]_m, $$
where  $B(\kappa)= L_{1}(\kappa)/L_{0}(\kappa), ~~L_{p}(\kappa)= \sum_{r=0}^{m-1}\cos\left(p\frac{2\pi r}{m}\right)e^{\kappa \cos\left(\frac{2\pi r}{m}\right)}, p\in \mathbb{Z}_m,$ and $[x]_m$ denotes  the closest integer to $x$, modulo $m$.

Asymptotic normality of these estimate follow using the results in  \citet{pewsey2004} for the case (a) and  \citet{mardiaetal2016ar} for the cases (b) and (c). 
  
{\bf Testing of Hypothesis.}

 We will  give more details in the next section on testing of hypothesis as we apply the methodology as required in the next section but we outline some main points.

{\it The Rayleigh Test}  is well known to test uniformity  under the von Mises distribution, that is  to test ( see, for example, \citealt{mardiajupp2000})
$$ H_0 : \kappa =0,~~ H_1: \kappa >0, $$
where the mean $\mu$ is unknown. It is based on $T_1^2= 2n \bar{R}^2$, which under $H_0$, is approximately chi-squared with 2 degrees of freedom.
Now, given the data vector $\boldsymbol{w}$ of iid observations  under CDVM, the log-likelihood ratio test statistic ($T$)  can be written as
$$ T(\boldsymbol{w},\hat{\kappa},\hat{\mu})= -2 LL(\boldsymbol{w}\vert m, \kappa=0)  + 2 LL(\boldsymbol{w}\vert \hat{\kappa},\hat{\mu}),
$$
where $(\hat{\kappa}, \hat{\mu})$ is the mle based on the data vector $\boldsymbol{w}$ and $LL(\cdot)$ is the log-likelihood. If we denote the computed value of $T$ in the data sample by $T_{d}$, then 
 $$ \mbox{p-value}=P\left(T(\boldsymbol{w},\hat{\kappa},\hat{\mu} ) \geq T_{d}\right).$$  The test rule is then to reject $H_0$ for values of the p-value in comparison to a chosen significance level and  
use bootstrap  for the mle and the p-values for the tests. 
 For large $m$, the T-test = the Rayleigh test. We can extend it  easily to the circular location family including CDWC  which is easier to use as the normalizing constant is simpler.

\section{Examples}
\label{S:examples}
We apply some of the models developed in the previous sections to analyze real data  on  roulette wheel outcomes and smart health monitoring readings on SBP acrophase. In both these situations, data are circular and discrete. Also, in both these situations, data are generated in abundance daily, although they may not usually be accessible in the public domain. Interestingly, the acrophase  is an example where the data has an irregular discrete support.  Our analysis and findings that are presented below are mainly illustrative of the kind of insights that are possible through the different circular discrete models.
\subsection{Roulette wheel data: online gaming and casino spins}
\label{S:examples_roulette}
There has been ongoing  search to find a plausible test for testing unbiasedness of a roulette wheel. The problem is now more pressing than ever before with the rise of  many online gaming sites, e.g. \url{ https://10bestcasinos.co.uk/en-en_d_rl.html }. For example, the UK Gambling Commission requires statistical testing to ensure fairness, by an approved third party as per the guidelines provided in their ``Testing strategy for compliance with remote gambling and software technical standards" at \url{http://www.gamblingcommission.gov.uk/}.)

Indeed, \citet{pearson1894,pearson1897a}  was captivated by this problem and had acquired his  data of   $n=16,563$ roulette spins from the Monte Carlo  Casino  as recorded in a journal
``La Monaco". He constructed three tests, but as the subject of Directional Statistics was still developing, it was usual to ignore the circular aspect of the data. A brief historical insight into his work, along with images of a typical European roulette wheel, are given in  the supplement. Karl Pearson's original sequence of roulette spin data is not available and  we work with data from an online roulette simulator as well as some industrial casino data obtained from spins of  different European roulette wheels. Note that here  the number of outcomes is  $m=37$ (as against an American roulette $m=38$). In some cases, the data is available as a streaming sequence  of outcomes from successive spins of the roulette wheel, and in others we may just have the frequency distribution of the outcomes without knowledge of the sequence. Accordingly, we illustrate different types of analysis. 
\subsubsection{Analysis of streaming sequence of roulette outcomes }
\label{S:examples_roulette1}
We look at three different sequential roulette data:
\begin{itemize}
\item[]{\bf Roulette data 1} Our first data is of size  $n=1000$, a sequence of outcomes from successive spins of an online European roulette simulator available at \\(\url{http://datagenetics.com/blog/july12015/index.html}) (last accessed 9-sept-2020). 
\item[] {\bf Roulette data 2} This  data has  outcomes from successive spins of a real European roulette recorded in a casino  in Slovenia. 
\item[] {\bf Roulette data 3} This data has  outcomes from successive spins  from the same casino as roulette data 2 , but from a different roulette wheel. 
\end{itemize}

Roulette data 1, 2 and 3 are of size $n=1000,~ 8299 $ and  $~8106$ respectively and 
rows (i)-(iii) of Table \ref{T:roulettefreq}  give their frequency distributions. The main challenge is to detect a possible bias in a roulette based on a streaming sequence of spin outcomes.  We note here that since the roulette data 1,2 and 3 are available as a time series, we carried out a  test for serial independence, in the lines of \citet{watson_beran1967}, but adapting to discrete data (more details are in \citet{mardiasriram2020ar}). The test indicates that overall there is no dependence.  Our analyses and findings for these data are as follows.

\begin{table}[htbp]
  \centering
  \caption{\label{T:roulettefreq} Frequency distributions for the online roulette data 1,2,3 and 4.  The entries of angular positions ($r$) in the first row  correspond to angles $2\pi r/m$ on the circle, the second row shows the corresponding label on the roulette wheel and the other rows show the frequencies of outcomes.}
  \resizebox{\textwidth}{!}{
    \begin{tabular}{cp{1.5cm}|ccccccccccccccccccc}
  ~~~  &  ~~$r$~~& ~~0 & ~~1 & ~~2 & ~~3 & ~~4 &~~5 &~~6 & ~~7 &~~8 & ~~9 & ~10 &~11 &~12 &~13 &~14 &~15 &~16 & ~17 & ~18\\
   \hline
    \end{tabular}}
  \resizebox{\textwidth}{!}{
   \begin{tabular}{cp{1.5cm}|gbababababababababa}
  ~~~  & label     & ~~0     & ~26    & ~~3     & ~35    & ~12    & ~28    & ~~7     & ~29    & ~18    & ~22    & ~~9     & ~31    & ~14    & ~20    & ~~1     & ~33    & ~16    & ~24    & ~~5 \\
    \hline
    \end{tabular}}
    \resizebox{\textwidth}{!}{
    \begin{tabular}{cp{1.7cm}|ccccccccccccccccccc}
(i)& Data 1  & 45    & 19    & 34    & 20    & 29    & 25    & 25    & 19    & 34    & 16    & 29    & 25    & 32    & 31    & 32    & 18    & 29    & 32    & 28 \\
  (ii)&  Data 2  & 242   & 228   & 231   & 213   & 241   & 211   & 202   & 221   & 208   & 216   & 230   & 199   & 232   & 230   & 220   & 208   & 196   & 215   & 203 \\
  (iii)& Data 3 & 230   & 202   & 184   & 231   & 200   & 221   & 169   & 187   & 224   & 212   & 232   & 221   & 255   & 236   & 225   & 190   & 189   & 199   & 214 \\
   (iv)  & Data 4 & 76 & 78 & 72 & 85 & 74 & 101 & 93 & 71 & 67 & 70 & 67 & 103 & 104 & 102 & 93 & 82 & 84 & 73 & 65   \\
   \hline
    \end{tabular}}
   ~\\
   \resizebox{\textwidth}{!}{
     \begin{tabular}{cp{1.5cm}|cccccccccccccccccc|c}
    ~~~& ~~$r$~~& ~19 & ~20 & ~21 & ~22 & ~23 &~24 &~25 & ~26 &~27 & ~28 &~29 &~30 &~31 &~32 &~33 &~34 &~35 & ~36~ & Total \\
    \hline
    \end{tabular}}
\resizebox{\textwidth}{!}{
    \begin{tabular}{cp{1.5cm}|bababababababababa|c}
 ~~~& label     & ~10    & ~23    & ~~8     & ~30    & ~11    & ~36    & ~13    & ~27    & ~~6     & ~34    & ~17    & ~25    & ~~2     & ~21    & ~~4     & ~19    & ~15    & ~32    &~~~~~~~\\
     \hline
    \end{tabular}}       
\resizebox{\textwidth}{!}{
    \begin{tabular}{cp{1.4cm}|cccccccccccccccccc|c}
  (i)&  Data 1  & 27    & 25    & 31    & 38    & 33    & 32    & 26    & 11    & 21    & 24    & 27    & 21    & 29    & 33    & 31    & 22    & 23    & 24    &1000\\
 (ii) &    Data 2 & 244   & 260   & 213   & 197   & 219   & 240   & 244   & 242   & 227   & 221   & 239   & 233   & 236   & 227   & 229   & 221   & 218   & 243 & ~8299\\
 (iii)&   Data 3& 196   & 203   & 206   & 207   & 213   & 247   & 222   & 199   & 201   & 233   & 250   & 254   & 274   & 239   & 220   & 256   & 218   & 247 & ~8106\\
     (iv)& Data 4 & 90 & 80 & 73 & 101 & 75 & 86 & 79 & 89 & 87 & 75 & 91 & 85 & 94 & 106 & 89 & 78 & 86 & 70 & 3094 \\ 
    \hline
    \end{tabular}}            
  \end{table}

 {\bf Analysis 1} : First, we carry out bias testing, which is equivalent to testing for uniformity, i.e.  $H_0: \tau=0$ (unbiased wheel) vs. $H_1:\tau \ne 0$ (biased wheel). For this purpose, we use a log-likelihood ratio  test statistic (denoted $T$), which is computed as the  difference in the log-likelihoods at the maximum likelihood estimates (mle) and at the null hypothesis $\tau=0$. Recall for the continuous circular location scale family, $\tau=0$ corresponds to the uniform distribution on the circle, so a specific value of $\mu$ is not required for $H_0$. We carry out this analysis using the CDWC model. Table \ref{T:MLE_test_PDWC} shows the mle, test statistic and the p-value for each the roulette data 1,2 and 3. Comparing the p-value with a 5\% significance level, we conclude that the evidence for bias does not exist for roulette data 1, is weak for roulette data 2 and strong for roulette data 3. The estimated mode for data 3 is $\hat{\mu}=5.34$ approximately corresponding to the angular position $r\approx 31$. See Supplement  for the circular histogram of data 3. It turns out these conclusions happen to be the same if we had used the CDVM model, or if we had used some alternative tests known in the context of continuous data (see Supplement  for details). 
\begin{table}[htbp]
  \centering
  \caption{\label{T:MLE_test_PDWC} Based on CDWC model, results of Analysis 1 for roulette data 1-3. } 
  
  \resizebox{.9\textwidth}{!}{  \begin{tabular}{cccccccccccc}
  & roulette     & n      & $\hat{\rho}$ & SE($\hat{\rho}$) & $\hat{t}=\lfloor \frac{m\hat{\mu}}{2\pi}\rfloor$& $\hat{\theta}=\hat{\mu}$ &$SE(\hat{\theta})$ & $T$     & $SE(T)$ & p-value \\
    \hline
      (i)   & data 1 & 1000    &         0.019  &         0.017  & 17    &         2.976  &        1.236  &         0.688  &         3.010  &      0.711 \\
    (ii)  & data 2 & 8299   &         0.020  &         0.007  & 30    &         5.106  &        0.414  &         6.676  &         5.420  &      0.046 \\
    (iii) & data 3 & 8106   &         0.030  &         0.008  & 31    &         5.340  &       0.279  &      14.333  &         7.924  &                0.000   \\
    \end{tabular}}  
\end{table}%

{\bf Analysis 2} : Analysis 1 does not use the streaming nature of the outcomes. Here, we delve further into the sequence of outcomes.
Let us denote the sequence of angular positions of the roulette outcomes by $\left\{w_i : i\in\{1,2,\ldots, n\}\right\}$. The mapping of angular positions to labels on the roulette wheel is given in Table \ref{T:roulettefreq}. 
 Our goal is to estimate a change-point in the data. The model with a changepoint at $i=K$ can be constructed as follows:
\begin{equation}
\label{eq:pmf_chngpt}
w_i \sim \begin{cases} p(\cdot \vert m, \tau_1, \mu_1) ~~~\mbox{ if } i\leq K \\  p(\cdot \vert m, \tau_2, \mu_2) ~~~\mbox{ if }i>K, \end{cases}
\end{equation}
where $p(\cdot\vert m, \tau, \mu)$ is the probability function of the discrete circular distribution as in Equations (\ref{eq:loc_fam_marg}) and  (\ref{eq:PDgen}). Of particular interest  is to detect a change from uniformity, i.e. $\tau_1=0$, where $\mu_1$ can be arbitrary but we take it as $0$ without loss of generality. The likelihood for the data can then be written as
\begin{equation}
\label{eq:lik_chngpt}
L(\boldsymbol{w}\vert (\tau_1=0,\tau_2), (\mu_1=0, \mu_2), K) = \prod_{i=1}^K  p(w_i \vert m, \tau_1=0, \mu_1=0) \times \prod_{i=K+1}^n  p(w_i \vert m, \tau_2, \mu_2).
\end{equation}
We can apply a Bayesian approach with standard Markov Chain Monte Carlo methods for estimation, using non-informative flat priors on the unknown parameters, viz $\tau_2,\mu_2$ and $K$.  For our data, we carry out the changepoint analysis using the $CDWC(m, \rho, \mu)$ model, first based on the full data sequence and then based on partial data sequences. Specifically, we use a Gibbs sampling procedure to obtain the posterior distributions,  by taking the support of $(\rho_2,~\mu_2,~K)$ to be a  suitably fine grid of values, viz. $\rho_2 \in \{0, 0.001, 0.002,\ldots, 0.999\}$, $\mu_2 \in  \{0, 0.001\times 2\pi, 0.002\times 2\pi,\ldots, 0.999\times 2\pi\}$ and $K\in \{1,2,\ldots, n\}$. The posterior distribution summaries obtained for each of the  full sequences of roulette data 1,2, 3 are shown in Table \ref{T:chngpt_ar04hispar0204}. For any given roulette data sequence, we conclude that there is evidence for a changepoint  if the 95\% highest probability density (hpd) credible interval for $\rho_2$ is removed from $0$. Accordingly, we conclude that there is no evidence for  a change from uniformity in data sequences 1 and 2. For roulette data 1 and 2, we also see that the 95\% hpd interval for $K$ spans a very large range between $1$ and $n$, as one would expect if the distribution of $K$ is close to a discrete uniform on $\{1,2, \ldots, n\}$. This further supports the absence of a changepoint in roulette data 1 and 2. 

However, there is evidence for a changepoint in roulette data 3 since the 95\% hpd interval for $\rho_2$ is clearly  removed from 0. We can also see that the estimated posterior mode for the changepoint in roulette data 3 is  $K= 1226$. Therefore, subsequent to the changepoint, the distribution of outcomes in roulette data 3 has changed from uniform distribution to one that has a single mode (at $\mu_2=5.305$ with angular position  $r\approx$ 31).  Such a bias might be resulting from a slight ``tilt" in the roulette wheel downwards to favour such a mode. In Section \ref{S:roulettefixedframe}, we also look at a different type of bias, possibly resulting from ``wobble" of a roulette wheel.
\begin{table}[htbp]
  \centering
  \caption{\label{T:chngpt_ar04hispar0204}: Posterior distribution summaries from the $CDWC$ changepoint model, for parameters $(\mu_2, \rho_2, K)$ based on full data ranges (i.e. $1:n$) for roulette data 1,2, 3. The distribution before changepoint is assumed to be uniform, i.e. $\rho_1=0, \mu_1=0$.   
  }
   \resizebox{\textwidth}{!}{
    \begin{tabular}{c|cccc|cc|cccc|c}
    \toprule
             &\multicolumn{4}{c|}{$\rho_2$}    & \multicolumn{2}{c|}{ $\mu_2$} & \multicolumn{4}{c|}{$K$} \\
             \hline
         Roulette & mean & sd    & \multicolumn{2}{c|}{95\% hpd interval}  & mean & sd & n& mode & sd    & \multicolumn{1}{c|}{95\% hpd interval}&changepoint \\
          \hline
          data 1 & 0.071 & 0.113 & [~~~~0~~~,& 0.288] &  2.747 &       1.822   & 1000&  998   & 305.34    &[~~91,~1000] & No \\
    
    data 2  & 0.031 & 0.060 & [~~~~0~~~,& 0.086] & 5.338& 0.949 & 8299& 8294   & 2828.61 & [244   ,~8298]   & No \\

data 3   & 0.035 & 0.010 & [0.015~~,& 0.053] & 5.305 & 0.270 &  8106&     1226  &      723.13  &  [~~~3,~2612]& Yes, at 1226\\
  
    \bottomrule
    \end{tabular}}%
\end{table}%

\begin{figure}[!htb]
    \begin{subfigure}{.5\textwidth}
 \center
  \caption{$\rho_2$}
	  \includegraphics[width=1.5in, height=3in, angle=270]{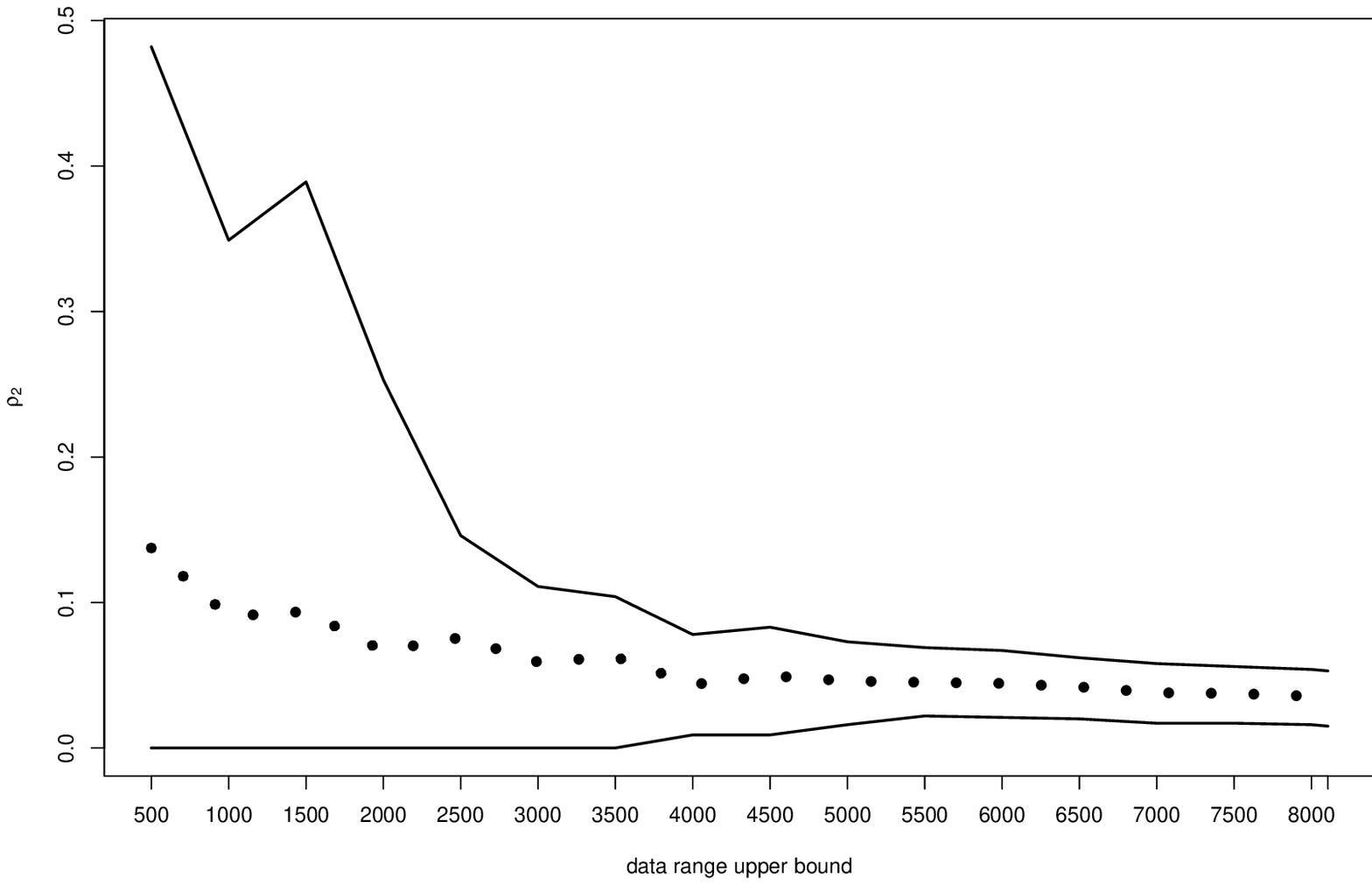}
	\end{subfigure}   
	  \begin{subfigure}{.5\textwidth}
 \center
  \caption{$K$}
	  \includegraphics[width=1.5in, height=3in, angle=270]{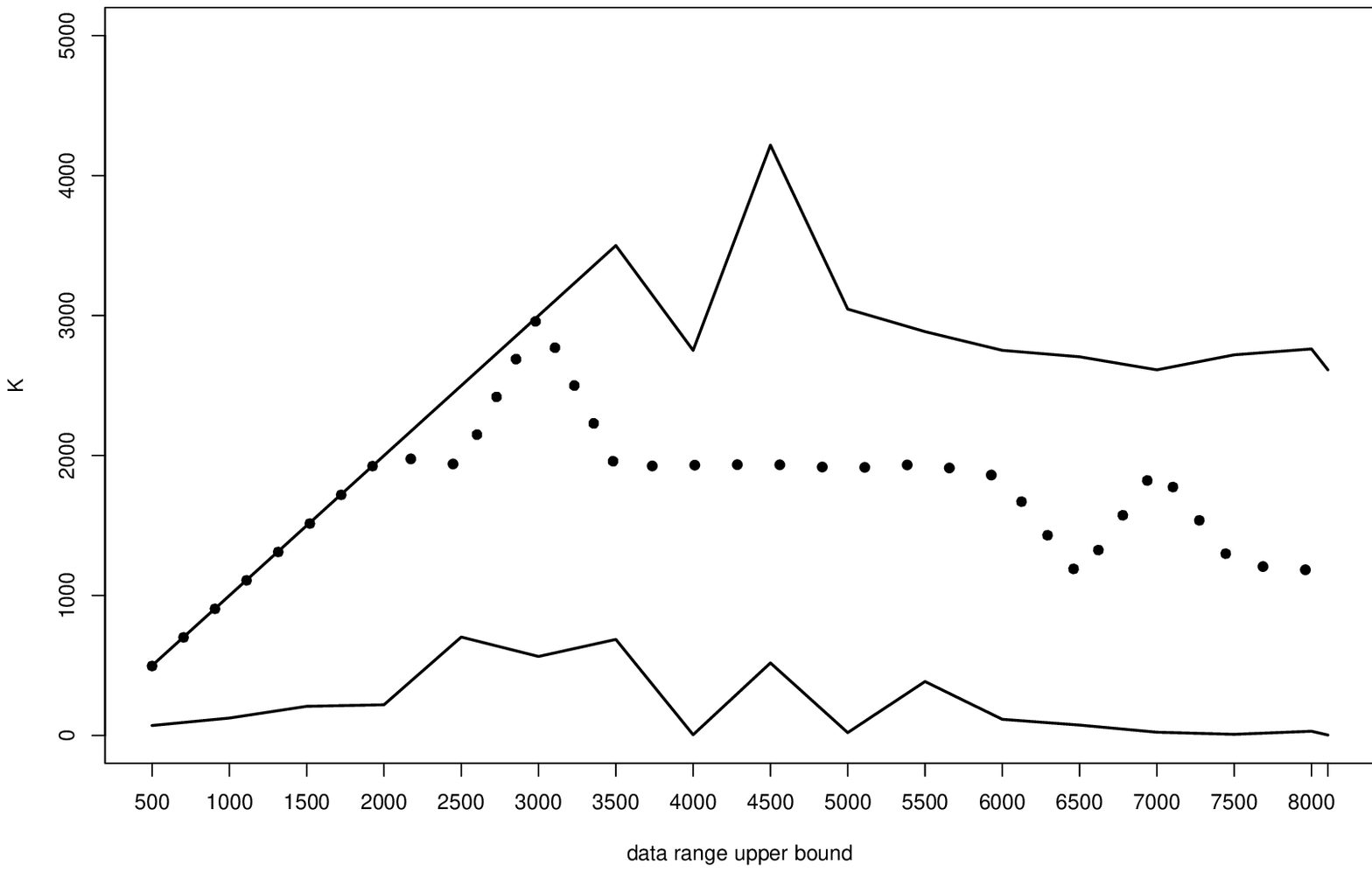}
	\end{subfigure}   
	 \caption{\label{F:out_roulette_chngpt_figdata3} Results of changepoint analysis based on different partial data ranges for Data  3. 
plot (i) shows posterior summaries for $\rho_2$ and (ii) for $K$ plotted against data range upper bound. In each plot, solid  lines mark the 95\% credible intervals. The  dotted line in plot (i) shows the posterior mean for $\rho_2$ and in plot (ii) shows the posterior mode for $K$.   }
          \end{figure}

{\bf Testing streaming consistency.} Further, to see how early the changepoint would have been detected, we treat the spins as streaming data by increasing the number of spins sequentially by 500. So, we  apply the same changepoint detection procedure on partial data sequences of roulette data 3, i.e. spins $1:u$ for  $u \in \{500, 1000, 1500, \ldots, 8000, 8106 \}$. Figure  \ref{F:out_roulette_chngpt_figdata3}  shows the plots of the 95\% hpd intervals, posterior mean for $\rho_2$ and mode for $K$  for different choices of the upper bound.  We can see that starting from an upper bound of $4000$ onwards, the 95\% credible interval for $\rho_2$ appears to be removed from 0, and for $5500$ onwards it is even more clearly removed from $0$. Correspondingly the posterior mode for $K$ is settled somewhere between $1000$ and $2000$. Also, the 95\% hpd for $K$ appears to be stabilize for the data range 1:5500 and after. So,  while we start detecting the change weakly based on the first 4500 outcomes,  the evidence becomes stronger as we start including outcomes 1:5500 and after. 
          
We now comment on how our work is basically different from some other work on circular change point detection. \citet{pewsey_eduardo2020} have given a survey of changepoint detection with continuous angular data, but our method is for discrete circular data. There has been other work in this area motivated by control charts, e.g. \citet{lombard_maxwell2012},  \citet{lahagupta2011} and \citet{girijarao2020}. Compared to these approaches, our strategy is  different as it is  model-based with unknown parameters, whereas the method in \citet{lombard_maxwell2012} is nonparametric. Further, control charts, e.g.   \cite{lahagupta2011} and \cite{girijarao2020}, are not suitable for testing uniformity for the following reasons. Broadly speaking, in their method, given a distribution, $n$ samples of size $n_1$ are generated. For each sample, based on the $n_1$ observations, the circular mean is computed. Then two quantities are found after sorting the  $n$ values, viz. (a) CCR (clockwise control ray) by eliminating the first $\alpha$ percent values and (b) ACR (anti-clockwise control ray) by eliminating the last $\alpha$ percent values. Such an approach may not be suitable for checking deviation from uniformity because the ACR and CCR will be wide and symmetric around $\pi$ and so it will not have any observations that fall in the rejection region if the data is actually coming from a distribution concentrated around $\pi$.  Although their work is also based on some model assumptions, it requires a priori fixing of parameter values to determine the ACR and CCR.
\subsection{Analysis of fixed frame of roulette outcomes}
\label{S:roulettefixedframe}
Here, we consider an example where only a fixed frame of roulette outcomes is available, i.e. instead of a streaming sequence we only have the frequency distribution of outcomes. As part of an industrial consulting project at University of Leeds, \cite{baines1990} recommended a protocol for certifying a casino roulette wheel as ``unbiased" using  five different statistical tests, two based on chi-square statistic and the other three based on variations of the Rayleigh test (see Supplement  for detailed description of these tests).  As per the protocol, if all five statistical tests resulted in the acceptance of the Null Hypothesis ($H_0$) (i.e. no evidence of bias), then the wheel would be certified as having passed the randomness test.  However, if any of the five statistical tests resulted in the acceptance of the Alternative Hypothesis ($H_1$)(i.e. evidence of bias), a second series of at least $370$ spins would be collected and the five statistical tests repeated on the new data as well as the combined data from the two series. The choice of the level for the tests would be contextual, i.e. chosen between 5\% and $0.1\%$ depending on the acceptable risk of wrongly concluding that the wheel is biased. The cost implications also influence the number of runs to be considered for such tests. 

Based on an analysis of different sub-series of roulette spins as well as the combined series, Al Baines reported evidence to suggest that the roulette had a quadri-modal bias with possible modes at the roulette slot positions 12,20,30 and 21,  which correspond to angular positions $r\in \{4,13,22,32\}$. However, his methodology is not based on a statistical model, so it does not systematically identify the nature of the specific alternative in $H_1$, in particular the positions of modes and their likelihood.

\begin{figure}[h]
\center
 \caption{\label{F:pdfs_fit} CDTS model fit to the combined series of roulette outcomes in \citet{baines1990} }
\includegraphics[width=2in, height=5in, angle=270]{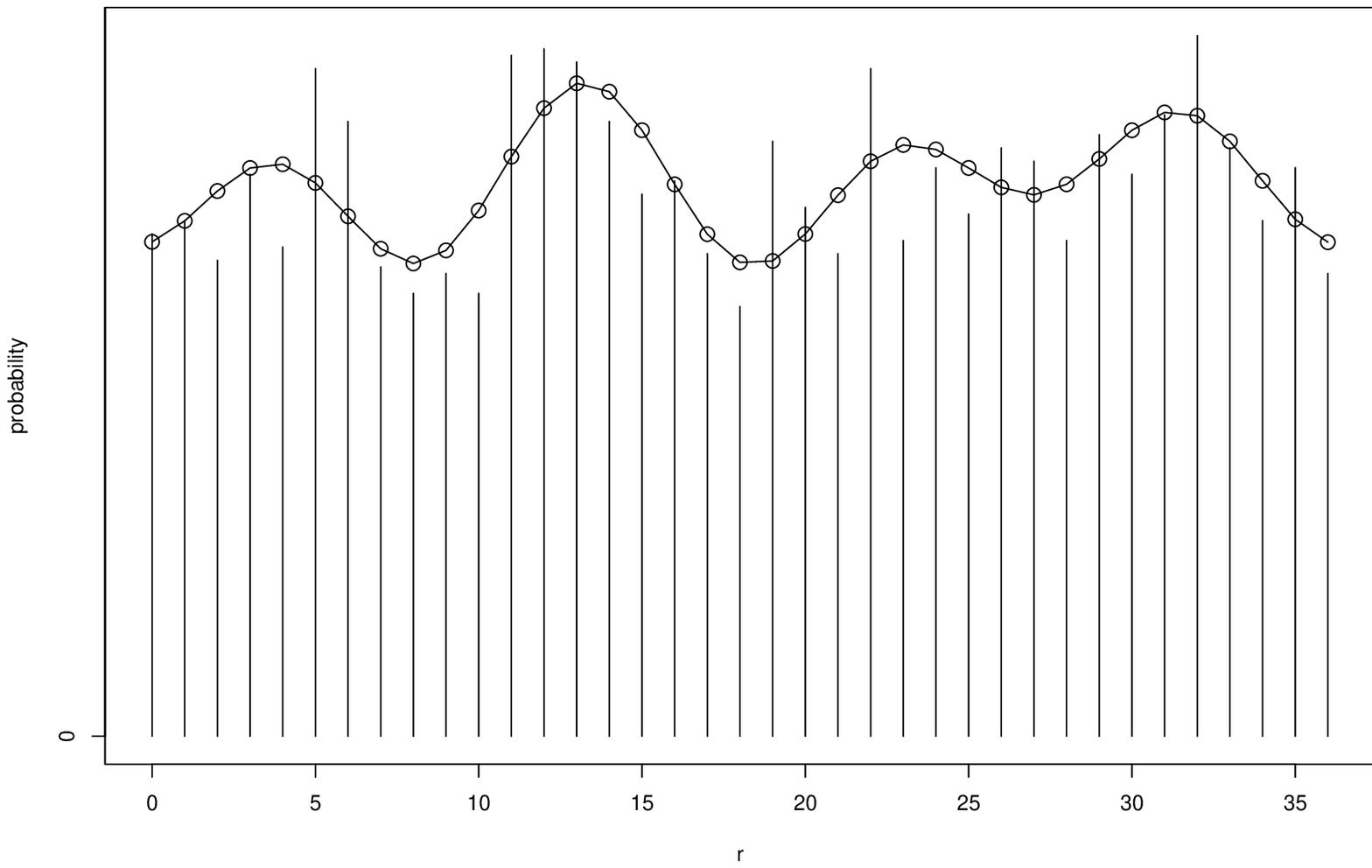}
\end{figure}

To allow some flexibility to capture the multimodal aspect of this data and to achieve the highest log pseudo marginal likelihood (see e.g. \citealt[p. 215]{carlin_louis2009}), we chose $J=4$ to fit the CDTS model as in Equation (\ref{eq:CDTS})  to the  combined series of roulette outcomes  used in \citet{baines1990}. We refer to this as `Data 4' and give its frequency distribution in row (iv) of Table \ref{T:roulettefreq}.
 We find it convenient to estimate the model using a Bayesian approach with a standard MCMC approach based on a random walk Metropolis-Hastings algorithm. Figure \ref{F:pdfs_fit} shows the fitted model along with the frequency distribution of the data.  The fitted model  suggests modes at angular positions $r\in \{4, 13, 23, 31 \}$, which is somewhat consistent with Al Baines' conclusion of 4 equi-spaced modes, although our finding suggests that the successive modes may not exactly be equi-spaced. As per our model, the angle subtended between successive modes (in degrees) are 87.6, 97.3, 77.8 and  97.3 respectively. The roulette slots corresponding to these angular positions  are 12, 20, 11, 2, with estimated probabilities 0.0279, 0.0319, 0.0289 and 0.0305 respectively.  So, our finding goes a bit beyond \citet{baines1990} to suggest that the roulette wheel possibly has an asymmetric wobble (i.e.  modes that are not exactly equi-spaced and with unequal probabilities ).  
\subsection{Acrophase data: ambulatory BP monitoring}
\label{S:acrophase}
\label{S:acrophase}
Systolic blood pressure (SBP) has a circadian rhythm. To monitor it, patients wear ambulatory devices that regularly measure and record blood pressure. In non-invasive smart health monitoring, these readings are typically recorded at predetermined, possibly irregular, discrete time points during daytime and nighttime. Our objective is to analyze ``acrophase", the time at which maximum SBP is attained in a day.  Monitoring the acrophase can provide an automated early warning of a possible
medical condition before it becomes clinically obvious. Typically, the readings taken by industrial monitoring devices are more frequent during daytime (e.g. each half hour) than at nighttime (e.g. each hour). For more details, see ~\\
\url{https://www.londoncardiovascularclinic.co.uk/cardiology-info/investigation/24-hr-ambulatory-blood-pressure }.~\\
Therefore, the resulting acrophase data is circular, discrete and supported on an irregular lattice on the circle. 

We use  acrophase data on SBP based on readings taken for an individual from 3-31-1998 to 7-7-2000, collected by the Halberg Chronobiology Center (University of Minnesota). We note that for a few days, there were multiple time points where maximum SBP was achieved. In such cases, we retained all such time points, thus resulting in a total of 880 data points. As mentioned in the introduction, typically, acrophase data is  extracted from SBP measurements at each half hour during daytime (8 am to 8 pm) and  each hour during nighttime (8 pm to 8 am).  If we map 8 am to 0 radians and 8 pm to $\pi$ radians, the acrophase times get mapped to an irregular support of $m=36$ points on the circle given by
\begin{eqnarray}
\left\{\frac{2\pi r}{48}, r =0,1,2,\ldots,24\right\} \bigcup \left\{\frac{2\pi r}{48}, r =26, 28,\ldots, 44, 46\right\},  \label{eq:supp_acro}
 \end{eqnarray}
Note that unlike the case of regular support where the points in the support are expressed in terms of $\mathbb{Z}_m$, here the points in the irregular support are expressed as $2\pi r/48$ to accommodate the half-hour and one-hour time points.
Table \ref{T:Kdata_freq} shows the frequency distribution of the data.  Our interest here is to estimate the centering and concentration parameters for this data.  While the centering parameter will be indicative  of the most likely timing of acrophase, the concentration parameter will indicate the extent of variability around that timing. We   adapt CD distributions  to an irregular support as follows and prefix such distributions   by   ICD, that is, for example, ICDVM stands for the CDVM with irregular support.  Let $\mathcal{S}=\{\theta_0, \theta_2, \ldots, \theta_{m-1}\}$ denote the  irregular circular lattice support (\ref{eq:supp_acro}). The conditionalized discrete   probability function for any  given parent pdf  $f(\cdot)$  on the irregular support  $\mathcal{S}$, ICD, is given by
\begin{eqnarray}
p(\theta_l)  =
\frac{f(\theta_l)}{\sum_{k=0}^{m-1} f(\theta_k)},  ~~~\theta_l \in \mathcal{S}. \label{eq:CDgen_irr}
\end{eqnarray}
  
\begin{table}[htbp]
  \centering
  \caption{\label{T:Kdata_freq} Frequency distribution of the  acrophase data. Also shown is $r$ corresponding to  each time point, which maps to the circle by $2\pi r/48$.}
    \resizebox{.85\textwidth}{!}{ 
       \begin{tabular}{lrrrrrrrrrrrr}
    Time & 08:00 & 08:30 & 09:00 & 09:30 & 10:00 & 10:30 & 11:00 & 11:30 & 12:00 & 12:30 & 13:00 & 13:30 \\
    $r$     & 0    & 1    & 2    & 3    & 4    & 5    & 6    & 7    & 8    & 9   & 10    & 11 \\
    freq  & 16    & 13    & 8     & 21    & 7     & 7     & 16    & 15    & 19    & 27    & 27    & 28 \\
    \end{tabular}%
    }
    ~\\~\\
        \resizebox{.85\textwidth}{!}{ 
        \begin{tabular}{lrrrrrrrrrrrr}
    Time  & 14:00 & 14:30 & 15:00 & 15:30 & 16:00 & 16:30 & 17:00 & 17:30 & 18:00 & 18:30 & 19:00 & 19:30 \\
    $r$    & 12   & 13    & 14    & 15    & 16    & 17    & 18    & 19    & 20    & 21    & 22    & 23 \\
    freq  & 31    & 27    & 27    & 30    & 24    & 34    & 46    & 58    & 54    & 76    & 79    & 50 \\
    \end{tabular}
    }
    ~\\~\\
      \resizebox{.85\textwidth}{!}{  \begin{tabular}{lrrrrrrrrrrrr}
    Time & 20:00 & 21:00 & 22:00 & 23:00 & 00:00 & 01:00 & 02:00 & 03:00 & 04:00 & 05:00 & 06:00 & 07:00 \\
    $r$     & 24    & 26    & 28    & 30    & 32     & 34     & 36     & 38     &40     & 42    & 44    & 46 \\
    freq  & 39    & 33    & 16    & 7     & 5     & 1     & 2     & 6     & 4     & 5     & 3     & 19 \\
    \end{tabular} }
    ~\\~\\
\end{table}%
\begin{table}[htbp]
  \centering
  \caption{\label{T:est_PDVMsimdata2} Estimated center and concentration parameters for data simulated from ICDVM with the same irregular support as acrophase data. The true values of parameters are  $\mu=0.785$, $\rho=0.6$ , $\kappa=1.516$.}
 \begin{subtable}{.45\textwidth}
  \center
  \caption{Sample statistics with bootstrap standard errors }
    \resizebox{\textwidth}{!}{
    \begin{tabular}{lccc}
     parameter  & estimate   & se    & 95\% interval \\
      $\mu$    &0.912	& 0.045 & [0.822, ~1.008] \\
      $\rho$ & 0.493	& 0.020	& [0.457,~0.532]\\
    $\kappa$ & 1.139 &	0.061& [1.031,~1.265] \\
    \end{tabular}}
    \end{subtable}
 \begin{subtable}{.45\textwidth}
  \center
  \caption{von Mises model  }
    \resizebox{\textwidth}{!}{
    \begin{tabular}{lccc}
     parameter  & estimate   & se    & 95\% interval \\
   $\mu$    & 	0.910 & 0.046 & [0.822, ~1.001]\\
$\rho$ & 0.491& 0.020& [0.451, ~0.529] \\  
   $\kappa$  &1.131 & 0.060 & [1.014,~1.253]\\
 \\
    \end{tabular}}
    \end{subtable}
    
 \begin{subtable}{.5\textwidth}
  \center
  \caption{ICDVM model on irregular support }
    \resizebox{\textwidth}{!}{
    \begin{tabular}{lccc}
     parameter  & estimate   & se    & 95\% interval \\
   $\mu$    & 0.812	& 0.033 & [0.747,~0.877]\\
$\rho$ & 0.608 &	0.015&	[0.578,~0.637] \\  
   $\kappa$   & 1.553 &	0.063& [1.429,~1.676]\\
 \\
    \end{tabular}}
    \end{subtable}

\end{table}%
\begin{table}[htbp]
  \centering
  \caption{\label{T:est_Kdata} Parameter estimates for ICDVM on irregular support for the acrophase data}
    \resizebox{.6\textwidth}{!}{
    \begin{tabular}{lccc}
     parameter  & estimate   & se    & 95\% interval \\
 $\mu$    & 2.462 ($\approx 17:30$ hrs) & 0.049 &[2.367, ~2.554] \\
    $\kappa$ & 1.114 & 0.060 & [0.994, ~1.230] \\
      $\rho$   & 0.485 & 0.020 & [0.444, ~ 0.522] \\
    \end{tabular}}
  \end{table}%

We first assess the effect of using VM vs ICDVM via simulations, mimicking the acrophase data,
that is the irregular support, and the same data size ($n=880$) and call it ``simulated acrophase data" where  we take the true value of ICDVM to be $\mu=0.785$ , $\kappa=1.516$, and $\rho=0.6$. These were selected to approximately match the centering and concentration parameter computed from the acrophase data.
Table \ref{T:est_PDVMsimdata2} shows the estimation of centering parameter ($\mu$) and concentration parameter ($\kappa$) for this simulated acrophase data. Parts (a) and (b) of  Table \ref{T:est_PDVMsimdata2} do not consider the discrete nature of the data. Part (a) uses sample statistics along  with bootstrap standard errors, namely the circular mean of the data for $\mu$ and the mean resultant length $\bar{R}$ for $\rho$. Part (b) estimates the parameters by assuming a continuous model, namely von Mises.   
  Both (a) and (b) are unable to  closely estimate the true $\mu$ and $\kappa$ ($\rho$) parameters unlike the discrete ICDVM model shown in part (c). Especially, the 95\% interval for $\kappa$ in parts (a) and (b) do not capture  the true value. Thus the inappropriateness of applying techniques that are otherwise meant for continuous data become even more apparent when we are dealing with discrete circular data on an irregular support. In summary, this simulation experiment   clearly illustrates that using methods meant for continuous data on discrete data can be misleading. 
 
 On the other hand, from Table \ref{T:est_Kdata} where the parameters are estimated using the acrophase data,  the approximate 99\% confidence interval for $\mu$ is $17:30$ hours $\pm 33$ minutes, slightly different from the mode of 19:00 hrs observed in the data in Table \ref{T:Kdata_freq}.
This can be the effect of skewness in the data, which has led us to explore the analysis with the parent as a skew circular distribution  since one of our objectives is  constructions that allow flexible families of plausible discrete circular distributions.
  \\

 We have selected  the  circular distribution of \citet{kato_jones2015}, which  has four parameters that control the first four trigonometric moments, leading to  unimodal symmetrical as well as skew distributions as particular cases. This family also has an analytically tractable normalizing constant and its pdf is given by
\begin{equation}
\label{eq:skew_wc1}
g_{KJ}(\theta)=\frac{1}{2\pi} \left( 1+ 2\gamma \frac{\cos(\theta-\mu)-\rho \cos\lambda}{1+\rho^2 -2\rho\cos(\theta-\mu-\lambda)}\right),~-\pi<\theta \leq \pi,
\end{equation}
where the parameters are constrained by  
\begin{equation}
\label{eq:constraints}
0\leq \rho <1, 0\leq \gamma \leq (1+\rho)/2, -\pi\leq \mu, \lambda \leq \pi, \mbox{ and  } \rho\gamma \cos\lambda \geq (\rho^2+2\gamma-1)/2.
\end{equation}

 The  conditionalized discrete distribution from (\ref{eq:skew_wc1}), is given by  the   probability function   
\begin{equation}
\label{eq:skew_PDWC}
p\left(r\vert m, \rho, \mu, \gamma, \lambda \right)= \frac{1}{D^\star} \left( 1+ 2\gamma \frac{\cos(\frac{2\pi r}{m}-\mu)-\rho \cos\lambda}{1+\rho^2 -2\rho\cos(\frac{2\pi r}{m}-\mu-\lambda)}\right), r\in \mathbb{Z}_m,
\end{equation}
where 
 \begin{equation}
 \label{eq:skew_Cs1}
 D^\star=m\left(1+2\gamma \rho^{m-1}\frac{ \cos(m(\mu + \lambda)-\lambda)- \rho^m\cos\lambda }{1+\rho^{2m}-2\rho^{m}\cos(m(\mu+\lambda)) }\right),
\end{equation}
with the same constraints on parameters as in (\ref{eq:constraints}). The normalizing constant $D^\star$ is derived in the supplement . We will call this family, the conditionalized discrete Kato-Jones  family, or briefly as CDKJ family. The CDWC is obtained as a special case when $\lambda=0$ and $\gamma=\rho.$
Note that the constraints ensure that the probability function in (\ref{eq:skew_wc1}) is positive and hence also for the discretized version (\ref{eq:skew_PDWC}). We can now obtain from equation (\ref{eq:CDgen_irr}), the probability function for ICDKJ.

Adapting the method of moments approach in \citet{kato_jones2015}, we can obtain the estimates to use in the probability function of ICDKJ.  
The moment  estimates are 
$$\hat{\mu}=2.248, \hat{\lambda}=0.816, \hat{\rho}=0.495, \hat{\gamma}=0.584$$

\begin{figure}  
\center
 \caption{\label{F:acrophase_freq} Histogram for the acrophase data along with the fitted ICDKJ probability function in solid line }
\includegraphics[width=3in, height=5in, angle=270]{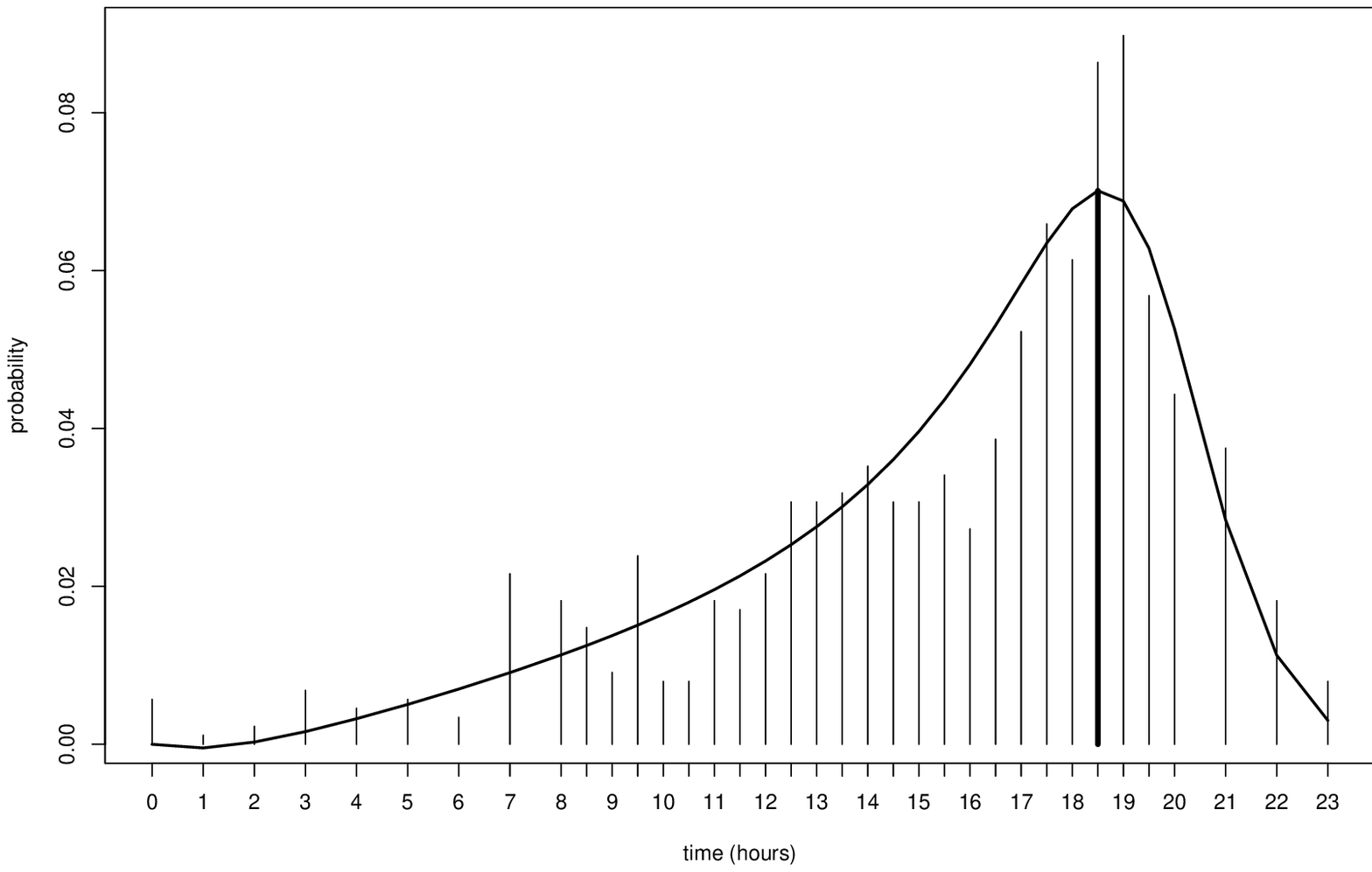}
\end{figure}
Figure \ref{F:acrophase_freq} shows the  histogram for the acrophase Data along with the fitted ICDKJ probability function with these estimates. It can be seen that 
the mode of the fitted ICDKJ distribution is 18:30 hours  approximately, which is roughly what is seen in the histogram. Also, visually,  the fit captures the skew behaviour   in the data adequately.

\section{Model Misspecification}
\label{S:comparisons}
In this section, we present some insights into the effect  of model misspecification. 
Suppose a true model is a given discrete model but we apply statistical methodology assuming that the data has come from a continuous circular distribution. To simplify the simulation experiments, we will assume that the continuous distribution is the parent distribution of a given  discrete distribution. We will study three situations:
\begin{itemize}
\item[] Case 1. The effect on the basic summary statistics.
We carry out a simulation study in the supplement, where we sample from the marginalized and conditionalized wrapped Cauchy distributions with $\mu=0,\rho=0.5$, and varying values of $m$. We  study   the first two trigonometric moments. Our conclusion is that discretization matters for $m \leq 20$.

\item[] Case 2. Behaviour of  the maximum likelihood estimates.
We study two sub-cases, 2a the regular discrete case, and 2b the irregular discrete case.

 In Case 2a, we use CDVM (CDWC) as the true  model and VM (WC) as the misspecified model  and calculate mle for $m=10, 20 $ with different values of the concentration
parameter.   We find that the mle using the misspecified  model leads to   biased estimates especially when there is a high concentration parameter in the true model (see the supplement ). 
 
For Case 2b, we have already given evidence  in Section \ref{S:acrophase}  that the effect on the mle  for a misspecified model can be more serious with irregular data. 

\item[] Case 3. Power of the  test of  uniformity.
For testing uniformity on highly dispersed data, as we have for our roulette data (Section \ref{S:examples_roulette1}), the  test based on a discrete distribution, CDVM or CDWC,  and the Rayleigh test based on a continuous distribution will lead to similar conclusions (see the supplement).

\end {itemize}

We now give some general remarks
 \begin{itemize}
\item[] Remark 1.
A key takeaway from the above comparison studies is that discreteness of the data cannot be ignored in general and such data may need to be modeled using discrete circular distributions.  At  least, these discrete distributions provide a bench mark to assess
any loss incurred in using continuous distribution.
\item[] Remark 2.  MD distributions are more relevant for the rounded circular data, whereas for any naturally circular discrete  data, CD distributions are more appropriate.
\item[] Remark 3. For $\lim m \to \infty $, the CD and MD distributions tend to their parent distribution as these are  the Riemannian  sums.  

\item[] Remark 4. Some other full scale  comparative studies  are given in \citet{mardiasriram2020ar}, including for example,  a comparison study of choices among conditionalized discrete distributions, based on   divergence measures such as Kullback-Leibler, $L_1$ and $L_2$. It is found that conditionalized discrete distributions resulting from von Mises and wrapped Cauchy can be very different from each other; so one cannot be easily approximated by the other family. In contrast,  the conditionalized discrete wrapped normal  and conditionalized discrete wrapped normal are very close to each other; so for practical purposes may be interchangeable in data analysis.

\end {itemize}

\section{Discussion} 
\label{S:conclusion}

 We have proposed flexible families of circular discrete distributions encompassing well established continuous circular  distributions, such as von Mises and wrapped Cauchy.  
Our analysis of model misspecification (Section \ref{S:comparisons}) highlights the importance of using discrete circular models for discrete data.
 We have selected the marginalized and conditionalized approaches for our analysis, but  other constructions such as the centered wrapped families  can be explored further. Also one can further explore the Beran family, in particular $\mathcal{B}_2$ and $\mathcal{B}_3$ given in Section \ref{S:maxentropy}). We have derived some insightful theorems interrelating the different methods of constructions. In particular, we have given an interesting  characterization that under some regularity conditions, marginalized  and conditionalized discrete distributions will be the same if and only if the parent circular distribution is uniform (leading to the discrete uniform).  The marginalized and conditionalized families of distributions have a significant potential for further development  beyond Directional Statistics. For example,  we have also proved a characterization on  the line that under some conditions, these two approaches lead to the same discrete distribution if and only if the parent  is the exponential distribution (see supplement )

We note that  some properties of the parent distributions such as unimodality and symmetry are inherited by the circular MD and CD. Also, maximum entropy characterization of  the von Mises distribution is inherited by CDVM. However not everything carries over, namely, the normalizing constant for the CDVM depends on both the parameters,  the Rayleigh test is no longer the likelihood ratio test for CDVM (see, supplement), and for CDWC the trigonometric  moments are not as simple.

It is worth noting that Karl Pearson recognized that the roulette wheel data goes beyond coin tossing data experiments and raises some difficult inference problems in assessing unbiasedness (see Supplement  for more details). Perhaps due to  the unavailability of adequate discrete circular models , it did
not have much impact at the time. Also, bias in a roulette wheel can be  due to a ``tilt" or can due to ``wobble" as we have seen in our examples. It could be said that wherever there is a wheel, there is  inherently a  natural circular discrete data, for example, a wheel used in  some TV shows- Wheel of Fortune and other shows , in industries  (bicycle wheel, umbrella,  and so on) .

The field is full of new challenges in statistical methodology, for example, the marginalized and conditionalized approaches are amenable to  extensions . We have outlined some   extensions in \cite{mardiasriram2020ar}, namely,  alternative approaches to allow for an irregular lattice support, and  extensions to higher manifolds such as the torus. However, it turns out that regular discretization on the sphere is not  straightforward (see supplement) , and there can be multiple ways of constructing conditionalized discrete distributions.
 
 Our overall recommendation from this paper is that ``If you have circular discrete data, you should start with a discrete model".
\FloatBarrier

\section*{Supplementary material}

Supplementary materials may be requested from the authors. More details and explanation will be available in the forthcoming monograph \citet{mardiasriram2020}.

\section*{Acknowledgments} 
We wish to thank Arthur Pewsey  for his help on some R queries related to his book \citet{pewseyetal2013}, to Shogo Kato for confirming a query, and to Colin Goodall, Peter Green, John Kent, Florian Pfaff and  John Wootton for helpful discussions.  We are deeply grateful to Mihael Perman for kindly providing us with the real data (roulette data 2 and 3) related to a casino in Slovenia, and Germaine Cornelissen-Guillaume for the real data on acrophase. We appreciate Neil Spencer  for making us aware of some other relevant work on roulette.     The first author would also like to thank the Leverhulme Trust for the Emeritus Fellowship. 

\bibliographystyle{rss}
\bibliography{ref_PDWC}

\end{document}